\documentclass[mnsc]{informs3_no_remark}

\OneAndAHalfSpacedXI
\usepackage{natbib}
\usepackage{algorithm}
\usepackage{algorithmic}
\usepackage{bm}
\usepackage{graphicx}
\usepackage{subfigure}
 \bibpunct[, ]{(}{)}{,}{a}{}{,}%
 \def\bibfont{\small}%
\def\tr{\mathop{\text{tr}}\kern.2ex}

\def\R{{\mathbb R}}
\def\E{{\mathbb E}}

\def\cS{{\mathcal S}}
\def\cC{{\mathcal C}}
\def\cA{{\mathcal A}}

\def\KL{{\text{KL}}}

\def\td{\text{d}}

\long\def\comment#1{}

\def\tr{\mathop{\text{Tr}}}

\def\cS{{\mathcal{S}}}
\def\cM{{\mathcal{M}}}

\TheoremsNumberedThrough     % Preferred (Theorem 1, Lemma 1, Theorem 2)
\ECRepeatTheorems

\EquationsNumberedThrough    % Default: (1), (2), ...
\MANUSCRIPTNO{MS-0001-1922.65}

\begin{document}
%%%%%%%%%%%%%%%%

% Outcomment only when entries are known. Otherwise leave as is and
%   default values will be used.
%\setcounter{page}{1}
%\VOLUME{00}%
%\NO{0}%
%\MONTH{Xxxxx}% (month or a similar seasonal id)
%\YEAR{0000}% e.g., 2005
%\FIRSTPAGE{000}%
%\LASTPAGE{000}%
%\SHORTYEAR{00}% shortened year (two-digit)
%\ISSUE{0000} %
%\LONGFIRSTPAGE{0001} %
%\DOI{10.1287/xxxx.0000.0000}%

% Author's names for the running heads
% Sample depending on the number of authors;
% \RUNAUTHOR{Jones}
% \RUNAUTHOR{Jones and Wilson}
% \RUNAUTHOR{Jones, Miller, and Wilson}
% \RUNAUTHOR{Jones et al.} % for four or more authors
% Enter authors following the given pattern:
%\RUNAUTHOR{}

% Title or shortened title suitable for running heads. Sample:
% \RUNTITLE{Bundling Information Goods of Decreasing Value}
% Enter the (shortened) title:
\RUNTITLE{A Primal-Dual Approach to CMDP}

% Full title. Sample:
% \TITLE{Bundling Information Goods of Decreasing Value}
% Enter the full title:
\TITLE{A Primal-Dual Approach to Constrained Markov Decision Processes} 
%with Applications to Inventory and Queueing Control}

% Block of authors and their affiliations starts here:
% NOTE: Authors with same affiliation, if the order of authors allows,
%   should be entered in ONE field, separated by a comma.
%   \EMAIL field can be repeated if more than one author
\ARTICLEAUTHORS{%
\AUTHOR{Yi Chen}
\AFF{Department of Industrial Engineering \& Management Sciences, Northwestern University, Evanston, IL, 60208 } %, \URL{}}
\AUTHOR{Jing Dong}
\AFF{Columbia Business School, New York City, NY, 12007 }
\AUTHOR{Zhaoran Wang}
\AFF{Department of Industrial Engineering \& Management Sciences, Northwestern University, Evanston, IL, 60208}
% Enter all authors
} % end of the block

\ABSTRACT{ In many operations management problems, we need to make decisions sequentially to minimize the cost while satisfying certain constraints. One modeling approach to study such problems is constrained Markov decision process (CMDP). When solving the CMDP to derive good operational policies, there are two key challenges: one is the prohibitively large state space and action space; the other is the hard-to-compute transition kernel. In this work, we develop a sampling-based primal-dual algorithm to solve CMDPs. Our approach alternatively applies regularized policy iteration to improve the policy and subgradient ascent to maintain the constraints. Under mild regularity conditions, we show that the algorithm converges at rate $O(\log(T)/\sqrt{T})$, where $T$ is the number of iterations. When the CMDP has a weakly coupled structure, our approach can substantially reduce the dimension of the problem through an embedded decomposition. We apply the algorithm to two important applications with weakly coupled structures: multi-product inventory management and multi-class queue scheduling, and show that it generates controls that outperform state-of-art heuristics. %We also demonstrate that CMDP can be viewed as a relaxation to weakly coupled Markov decision processes.
}%

% Sample
%\KEYWORDS{deterministic inventory theory; infinite linear programming duality;
%  existence of optimal policies; semi-Markov decision process; cyclic schedule}

% Fill in data. If unknown, outcomment the field
%\HISTORY{}

\maketitle
%%%%%%%%%%%%%%%%%%%%%%%%%%%%%%%%%%%%%%%%%%%%%%%%%%%%%%%%%%%%%%%%%%%%%%

% Samples of sectioning (and labeling) in MNSC
% NOTE: (1) \section and \subsection do NOT end with a period
%       (2) \subsubsection and lower need end punctuation
%       (3) capitalization is as shown (title style).
%
%\section{Introduction.}\label{intro} %%1.
%\subsection{Duality and the Classical EOQ Problem.}\label{class-EOQ} %% 1.1.
%\subsection{Outline.}\label{outline1} %% 1.2.
%\subsubsection{Cyclic Schedules for the General Deterministic SMDP.}
%  \label{cyclic-schedules} %% 1.2.1
%\section{Problem Description.}\label{problemdescription} %% 2.

% Text of your paper here

\section{Introduction}

In many sequential decision-making problems, a single utility might not suffice to describe the real objectives faced by the decision-makers.
A natural approach to study such problems is to optimize one objective while putting constraints on the others. 
In this context, the constrained Markov decision process (CMDP) has become an important modeling tool for  
sequential multi-objective decision-making problems under uncertainty. 
A CMDP aims to minimize one type of cost while keeping the other costs below certain thresholds. 
It has been successfully applied to analyze various important applications, including admission control and routing in telecommunication networks,
scheduling for hospital admissions, and maintenance scheduling for infrastructures \citep{alt}. 
Due to the complicated system dynamics and the scale of the problem, exact optimal solutions to CMDPs can rarely be derived.
Instead, numerical approximations become the main workhorse to study CMDPs.
In this paper, we propose a sampling-based primal-dual algorithm that can efficiently solve a wide range of CMDPs. %to derive good operational policies.
 
One basic approach to solve the CMDP is to use a linear programming (LP) formulation based on the occupancy measure. This approach faces two key challenges in implementations: it requires knowledge of the transition kernel of the underlying dynamical system explicitly; it does not scale well as the state space and action space get large. 
An alternative approach is to apply the Lagrangian duality. In particular, by dualizing the constraints and utilizing strong duality, 
we can translate the CMDP into a max-min problem, where for a  given Lagrangian multiplier, the inner minimization problem is just a standard Markov decision process (MDP).
This approach allows us to solve the inner problem using standard dynamic programming based methods. 
It does not require direct knowledge of the transition kernel as long as we can estimate  the value functions from simulated or empirical data. 
In implementations, one would iteratively update the MDP policy and the Langrangian multiplier. 
The current development of this approach requires solving the MDP to get the optimal policy for each updated Langrangian multiplier (see, for example, \cite{bpl,rlcc}), which can be computationally costly. A more natural idea is to solve the MDP only approximately at each iteration. 
In this paper, we investigate this idea and show that at each iteration, we only need to do one iteration of policy update to achieve the optimal convergence rate (in terms of the number of primal-dual iterations). 
Compared to the existing algorithms utilizing Lagrangian duality, our primal-dual algorithm can be run at a much lower cost at each iteration.
We also demonstrate that our algorithm can be easily combined with many other approximate dynamic programming techniques, such as Monte Carlo policy evaluation,
TD-learning, and value function approximations \citep{Sutton}.  

A key ingredient of our algorithm is regularized policy iteration. The standard policy iteration includes two steps: policy evaluation and policy improvement. The policy evaluation step calculates the action-value function under a given policy. Then, the policy improvement step defines a new policy by taking the action that minimizes the action-value function. Through a Kullback-Leibler (KL) regularization term, the regularized policy iteration modifies the policy improvement step by reweighing the probability of taking each action via a softmax transformation of the action-value function. 
This modification allows us to view the policy update step as running mirror descent for the objective function in the policy space \citep{md}.
In addition, we update the Lagrangian multiplier using subgradient ascent, which also belongs to the family of mirror descent methods. This unified viewpoint makes the improved primal-dual algorithm possible.  Noticeably, many recent developments in reinforcement learning also benefit from regularization, which has been shown to improve exploration and robustness. For example, Trust Region Policy Optimization and Proximal Policy Optimization use KL divergence between two consecutive policies as a penalty in policy improvement \citep{trpo, ppo}. Soft-Q-learning uses Shannon entropy as a penalty in value iteration \citep{softq}. 
\cite{rmdp} propose a unified framework to analyze the above algorithms via regularized Bellman operator (see also \cite{Liu,Shani,wang2019neural} for convergence analysis of regularized policy iteration). 

In terms of applications of the algorithm, we study an important class of CMDPs which we refer to as weakly coupled CMDPs \citep{Singh}.
A weakly coupled CMDP comprises multiple sub-problems that are independent except for a collection of coupling constraints.
Due to the linking constraints, the scale of problem grows exponentially in the number of sub-problems. Hence, even in the case where each sub-problem is computationally tractable, it can be computationally prohibitive to solve the joint problem. Our primal-dual algorithm naturally helps break the curse of dimensionality in this case. In particular, the weakly coupled CMDP can be decomposed into independent sub-problems in the policy iteration step. In this case, the complexity only grows linearly with the number of sub-problems. We also comment that the weakly coupled CMDP can be viewed as a Lagrangian relaxation of the weakly coupled MDP \citep{Adelman}. 
Even though there is a relaxation gap between the two, as we will demonstrate in our numerical experiments, the (modified) 
policy obtained via CMDP can perform very well for the original MDP problem in the applications we considered. 
%{\color{red} [remark: We do not consider weakly coupled MDP in inventory problem and do not check the gap as well. ]} }

We apply the primal-dual algorithm to solve two classical operations management problems: inventory planning and queue scheduling. For the inventory planning problem, we consider a multi-product newsvendor problem with budget constraints \citep{MPNP}. We formulate this problem as a weakly coupled CMDP and study a small-scale instance where we can numerically solve for  the optimal policy.
We show that our policy can indeed achieve ${O}(\log(T)/\sqrt{T})$ convergence in this case, where $T$ is the number of iterations. 
For the queue scheduling problem, we consider a multi-class multi-pool parallel-server system
where the decision-maker needs to route different classes of customers to different pools of servers in order to minimize the performance cost (holding cost plus overflow cost). 
We allow the service rates to be both customer-class and server-pool dependent.
Since each pool only has a finite number of servers, the routing policy needs to satisfy the capacity constraints.
This optimal scheduling problem can be formulated as a weakly coupled MDP. We consider instances where it is prohibitive to solve for the optimal policy.
Applying the Lagrangian relaxation, we solve the resulting weakly coupled CMDP by combining our primal-dual algorithm with value function approximation techniques. We show that our method generates comparable or even better policies than the state-of-art policies.

\subsection{Literature review}
In this section, we review some of the existing methods/results to solve CMDPs. The goal is to clearly state the contribution of our work.
Most existing algorithms for CMDPs is adapted from methods for MDPs, and can be roughly divided into three categories: LP based approaches, dynamic programming based approaches (including policy iteration and value iteration), and the policy gradient methods.

One LP based approach utilizes the occupation measure, which is the weighted proportion of time the system spends at each state-action pair. 
The objective and constraints can be written as the inner products of instantaneous cost functions and occupation measure.  
The other LP based approach utilizes the dynamic programming principle and treats the value function (defined on state space) as the decision variables. 
In particular, the optimal value function of an MDP is the largest super-harmonic function that satisfies certain linear constraints determined by transition dynamics. 
For the CMDP, we obtain an LP by combing the dynamic programming principle with the Lagrangian duality.
These two LP formulations are dual of each other \citep{alt}. 
Many works on LP based approaches aim to find efficient ways to solve the LPs by exploiting the structures of some specific MDPs/CMDPs. 
For example, \cite{BO} use the ellipsoid method to derive efficient algorithms for problems with side constraints, including the traveling salesman problem and the vehicle routing problem. \cite{OFP} studies a linear fractional programming method to solve CMDPs.
More recently, \cite{Caramanis} propose two algorithms based on the column generation and the generalized experts framework, respectively.
These algorithms are shown to be as efficient as value iteration. Another challenge of LP based approaches is that we need to know the transition kernel up front and in explicit forms. Some recent developments aim to overcome this challenge.  For example, \cite{SPD} reformulate the LP of an MDP as a saddle point problem and use stochastic approximation to solve it. \cite{SPD} combine the saddle point formulation of LP with value function approximation and develop a proximal stochastic mirror descent method to solve it. However, most existing developments in this line focus on MDPs. 

The policy/value iteration stems from the Bellman operator on the value function, and converges to the optimal value function linearly. 
In implementations, the value function can be estimated via simulation or data. Thus this class of methods does not require knowledge of the transition kernel up front. 
For example, in reinforcement learning, we learn the transition dynamics of the system while solving for the optimal policy \citep{Sutton}.
There are many works that apply policy/value iteration to solve CMDPs.
For example, \cite{Gattami} formulates the CMDP as a zero-sum game and uses primal-dual $Q$-learning to solve the game. 
It proves the almost sure convergence of the algorithm, but does not establish the rate of convergence.  \cite{bpl} study CMDPs in the offline learning setting, 
and combine various dynamic programming techniques with Lagrangian duality to solve it.  
\cite{rlcc} extend the constraints of CMDPs to nonlinear forms and propose to solve it via Lagrangian duality and $Q$-learning as well. 
An $O(1/\sqrt{T})$ rate of convergence is obtained in both \cite{bpl} and \cite{rlcc}. 
However, their algorithms require solving for the optimal policy at each updated Lagrangian multiplier.
%which would require several iterations of policy/value updates. 
Our method can be viewed as an improved version of their algorithms. In particular, our algorithm only requires one policy iteration for each updated Lagrangian multiplier while achieving the same convergence rate. We comment that \cite{bpl} and \cite{rlcc} consider more complicated settings than the classical CMDPs studied in this paper. 
It would be interesting to extend our primal-dual algorithm to solve the more general  problems. %as in \cite{bpl,rlcc}.
%However, it is possible to apply our primal-dual algorithm in their settings, which achieves more efficiency. }

Both LP based methods and policy/value iteration can suffer from the curse of dimensionality when dealing with a large action space. 
The policy gradient based approaches alleviate the dimensionality issue by approximating the policy using a parametric family of functions.  
In this case, searching over the policy space reduces to a finite dimensional optimization problem. 
%Commonly used policy gradient algorithms include: REINFORCE and the actor-critic algorithm \citep{Sutton}. 
Several works combine this idea with Lagrangian duality to solve large-scale CMDPs. For example, \cite{Borkar} and \cite{Bhatnagar} combine actor-critic algorithms with policy function approximations. \cite{rcpo} use two-timescale stochastic approximation. Beyond duality, several other policy gradient based approaches to solve CMDPs have been developed. For example, \cite{cpo} propose a trust region method that focuses on safe exploration. \cite{IPO} develop an interior point method with logarithmic barrier functions. \cite{las} propose to use Lyapunov functions to handle constraints. However, the key challenge of  policy gradient based methods to solve CMDPs is that the corresponding optimization problems are non-convex.  In most cases, only convergence to a local minimum can be guaranteed and the convergence rates are often hard to establish. 

\subsection{Organization of the paper and notations}
The paper is organized as follows. We first introduce the CMDP and review some classical results that are relevant to our subsequent development in Section \ref{background}. 
We then introduce our algorithm in Section \ref{algorithm}, and show that the algorithm achieves the optimal convergence rate in Section \ref{theory}. 
In Section \ref{WCMDP}, we discuss how our algorithm can be applied to (approximately) solve weakly coupled CMDPs and weakly coupled MDPs. 
We then implement our algorithm to solve an inventory planing problem and a queue scheduling problem in Sections \ref{sec:inventory} and \ref{sec:queue} respectively. Lastly, we conclude the paper and discuss some interesting future directions in Section \ref{conclusion}. 

The following notations are used throughout the paper. For a positive integer  $K$, we denote $ [K]$ as the set $ \{ 1,2,\ldots,K \} $. For a vector $\lambda \in \mathbb{R}^K$, $[\lambda]_k$ denotes its $k$-th coordinate and $\| \lambda \| = (\sum_{k=1}^K [\lambda]^2_k)^{1/2}$ denotes its $L_2$ norm. Given two vectors $a, b \in  \mathbb{R}^K$, we say $a \le b $ if the inequality holds for each coordinate, i.e., $[a]_k \le [b]_k\ \forall k\in[K]$. Given a vector $x \in \R^K$, $[x]^+=(\max\{[x]_1,0 \},\ldots, \max\{[x]_K,0 \})$. Finally, given two sequences of real numbers $\{a_n\}_{n\ge 1}$ and $\{b_n\}_{n\ge 1}$, we say $b_n=O(a_n)$,  $b_n=\Omega(a_n)$, and $b_n=\Theta(a_n)$ if there exist some constants $C,C'>0$ such that $ b_n \le C a_n$, $ b_n \ge C' a_n$, and $C'a_n\le  b_n \le C a_n$, respectively. We also introduce the $\tilde{O}(\cdot)$ notation when we ignore the logarithmic factors. For example, if $b_n \le C a_n \cdot \log(n)$, we denote it by  $b_n=\tilde{O}(a_n)$.

\section{Constrained Markov Decision Process} \label{background}
%In this section, we start by providing a brief overview of CMDP. 
%We then introduce two classic methods to solve CMDP.

%Then in Section \ref{classic_method}, we introduce two LP-based methods to CMDP: occupation measure method and Lagrangian multiplier method. which are dual of each other. Both of them requires the transition probabilities explicitly and are not scalable. To overcome the difficulties, starting from the Lagrangian duality, we develop a novel sampling-based primal-dual algorithm in Section \ref{algorithm}. Our algorithmic framework can also be combined with other advanced ADP/RL or sampling techniques easily in order to solve large-scale problems efficiently. 
 
%\subsection{Constrained Markov Decision Process} \label{background}
We start by considering a discrete-time MDP characterized by the tuple $(\cS,\cA,P,\gamma,\mu_0)$. Here,
$\cS$ and $\cA$ denote the state and action spaces; $P =\{ P(\cdot|s,a)\}_{(s,a)\in \cS \times \cA}$ is the collection of probability measures indexed by the state-action pair $(s,a)$. For each $(s,a)$, $P(\cdot|s,a)$ characterizes the one-step transition probability of the Markov chain conditioning on being in state $s$ and taking action $a$. Function $c=\{c(s,a)\}_{(s,a)\in \cS \times \cA}$ is the expected instantaneous cost where $c(s,a)$ is the cost incurred by taking action $a$ at state $s$. 
Lastly, $\gamma \in (0,1)$ and $\mu_0=\{\mu_0(s)\}_{s\in \cS}$ are the discount rate and the distribution of the initial state, respectively. 
Given an MDP $(\cS,\cA,P,\gamma,\mu_0)$, a policy $\pi$ determines what action to take at each state. 
We define the expected cumulative discounted cost with initial state $s_0$ under policy $\pi$ as 
\begin{equation} \label{def:value_func}
V^{\pi}(s_0)=(1-\gamma)\cdot\E^\pi \big[\sum_{t=0}^\infty \gamma^t \cdot c(s_t,a_t) \big| s_0 \big],
\end{equation}
where $ s_t,a_t$ are the state and action at time $t$ and $\E^\pi$ denotes the expectation with respect to the transition dynamics determined by policy $\pi$. We further weight the expected costs according to the initial state distribution and define
\begin{align} \label{loss_def}
    C(\pi)=\E_{s_0\sim \mu_0} \big[V^{\pi}(s_0) \big].
\end{align}
Our goal is to minimize the cost $C(\pi) $ over a suitably defined class of policies. 

As an extension to MDP, the CMDP model optimizes one objective while keeping others satisfying certain constraints. Specifically, in addition to the original cost $c$, we introduce $K$ auxiliary instantaneous costs $d_k=\{d_k(s,a)\}_{(s,a)\in\cS \times \cA},\forall\ k \in [K]$. The goal of a CMDP is to find the policy that minimizes the cost defined in \eqref{loss_def} while keeping the following constraints satisfied
\begin{align} \label{constraint_def}
  D_{k}(\pi)=(1-\gamma)\cdot \E_{s_0\sim \mu_0} \Big[ \E^\pi \big[\sum_{t=0}^\infty \gamma^t \cdot d_k(s_t,a_t) \big| s_0 \big]   \Big] \le q_k,\ \forall\ k \in [K]. 
\end{align}
In order to make the expression more concise,
we define 
$D(\pi):=(D_1(\pi),\ldots,D_K(\pi))^\top$, $ q:=(q_1,\ldots,q_K)^\top$, and write the constraints in \eqref{constraint_def}  as $D(\pi) \le q$. 

We remark that the CMDP is only one modeling choice to model problems with multiple objectives/constraints.
This particular modeling choice turns out to enjoy a lot of analytical and computational tractability as we will discuss next.
CMDP is also closely connected to an important class of MDPs -- weakly coupled MDP. 
In particular, CMDPs can be viewed as a relaxation of weakly coupled MDPs \citep{Adelman}. We will provide more discussions about this in Section \ref{WCMDP}.

%We conclude this section with some comments on the policy space, which is critical to understand the main results of this paper. 

\subsection{Policy Spaces}
Solving a CMDP requires finding the optimal policy over a properly defined policy space, which is a function space. 
Imposing suitable regularity conditions on the policy space will facilitate the development of algorithms.
We next introduce a few classes of commonly used policies.
It is natural to require that all policies are non-anticipative, which means that the decision-maker does not have access to future information. 
Define the history at time $t$ to be the sequence of previous states and actions as well as the current state, i.e.,
$h_t:=(s_0,a_0,\ldots,a_{t-1},s_t)$. Then a non-anticipative policy can be viewed as a mapping from $h_t$ and $t$ to the action space.
We refer to such a policy as a \textsl{``behavior policy''}. 
If a policy only depends on the current state $s_t$ and time $t$ instead of the whole history $h_t$, 
it is called a \textsl{``Markov policy''}. 
For a Markov policy, if it is independent of the time index $t$, it is referred to as a \textsl{``stationary policy''}. When a stationary policy is a deterministic mapping from the state space to the action space, it becomes a \textsl{``stationary deterministic policy''}. We use $\Pi$, $\Pi_M$, $\Pi_S$, $\Pi_D$ to denote the space of behavior, Markov, stationary, and stationary deterministic policies, respectively. 

Given an arbitrary policy space $U$, we can further generate a new type of policy called \textsl{``mixing policies''} via an initial randomization. Specifically, let $\rho$ be a probability measure on $U$. 
Under a mixing policy on $U$ with mixing probability $\rho$, we first draw a policy, say $\pi_g$, from $U$ following the distribution $\rho$. 
Then $\pi_g$ is executed for $t=0,1,2,\cdots$. We denote by $\cM(U)$ the space of mixing policies constructed from $U$. An important special case is $\cM(\Pi_S)$, i.e., the space of mixing stationary policies. 
When allowing mixing operation, we incorporate the randomness of the initial mixing into the calculation of the accumulated cost. 
In particular, for $\pi \in \cM(U)$ with initial randomization $\rho$, 
\begin{align*} 
C(\pi)=(1-\gamma)\cdot \E_{s_0\sim \mu_0} \Big[ \int_U \E^{\pi_g} \big[\sum_{t=0}^\infty \gamma^t \cdot c(s_t,a_t) \big| s_0 \big] \rho(\td g)  \Big].
\end{align*}
%In studies of MDP, we typically do not consider the mixing policy space because it is not necessary. However, for CMDP, introducing mixing policy space simplifies the development of algorithms and theory significantly without causing essential difference, especially for the primal dual methods. We explain more about it later.Now we consider the relationships between various policy spaces. 

By definition, we note that 
$$
\Pi \supseteq \Pi_M \supseteq \Pi_S \supseteq \Pi_D,\  \mathcal{M}(\Pi_S) \supseteq \Pi_S.
$$
A class of policies $U$ is called  a \textsl{``dominating class''} for a CMDP, if for any policy $\pi \in \Pi$, there exists a policy $ \bar \pi \in U$ such that 
$$
C(\bar \pi)\le C( \pi),\ D_k(\bar \pi)\le D_k(\pi),\ \forall\ k\in[K].
$$
For CMDPs, when the instantaneous costs $c(\cdot,\cdot)$ and $d_k(\cdot,\cdot)$ are uniformly bounded from below,  $ \Pi_S$ is dominating \citep{alt}. 
The class of mixing stationary polices $\mathcal{M}(\Pi_S)$ is also dominating in this case (Theorem 8.4 in \citep{alt}).

\subsection{Classical Approaches to Solve CMDPs} \label{classic_method}
There are two classical approaches to CMDPs. We use CMDPs with finite state and action spaces as examples. 
The first method utilizes the occupation measure. Given a policy $\pi$, the occupation measure is defined as
\begin{align} \label{occupation}
\nu^\pi(s,a):=(1-\gamma)\cdot \E_{s_0 \sim \mu_0} \Big[ \sum_{t=0}^\infty \gamma^t P^\pi(s_t=s,a_t=a | s_0) \Big],    
\end{align}
where $P^{\pi}(\cdot,\cdot|s_0)$ denotes the probability measure induced by policy $\pi$ with initial state $s_0$. 
Note that the occupation measure is the weighted long-run proportion of time that the system spends at each state-action pair. 
%Similarly, we define the state occupation measure by sum $\nu_{\mu_0}^\pi(s,a) $ over the actions, i.e. 
%\begin{align} \label{state_occupation}
%d^{\pi}_{\mu_0}(s)=(1-\gamma)\cdot \E_{s_0 \sim \mu_0} \Big[ \sum_{t=0}^\infty \gamma^t P^\pi(s_t=s| s_0) \Big].
%\end{align} 
We can then express the accumulated costs in \eqref{loss_def} and \eqref{constraint_def} as  
\begin{align*}
 C(\pi)&=\sum_{(s,a)\in \cS \times \cA} c(s,a)\cdot \nu^\pi(s,a)\\
 D_k(\pi)&=\sum_{(s,a)\in \cS \times \cA} d_k(s,a)\cdot \nu^\pi(s,a), ~ \forall\ k \in[K].
\end{align*} 
Let $\mathcal{Q}$ denote the set of feasible occupation measures, i.e., for any occupancy measure $\nu\in \mathcal{Q}$  there exists a policy $\pi$ that leads to $\nu$. By Theorem 3.2 in \cite{alt}, $\mathcal{Q}$ can be represented by the collection of vectors $\{\nu(s,a)\}_{(s,a)\in \cS \times \cA} $ that satisfies the following system of linear equations:
\begin{align*} 
 \sum_{(s,a)\in \cS \times \cA} \nu(s,a) \Big(\text{1}(s=s')-\gamma P(s'|s,a)\Big)& =(1-\gamma)\cdot \mu_0(s'),\ \forall\ s'\in\cS,
 & \nu(s,a)\ge 0,\ \forall\ (s,a)\in\cS \times \cA,
\end{align*}
where $\text{1}(\cdot)$ is the indicator function. Then we obtain the following LP formulation of CMDP 
\begin{align} \label{LP}
  &\min \sum_{(s,a)\in \cS \times \cA} c(s,a)\cdot \nu(s,a) \\
  &\ \text{s.t.} \quad  \big\{ \nu(s,a)\big\}_{(s,a) \in \cS \times \cA }  \in  \mathcal{Q}, \notag \\
  & \quad \ \ \sum_{(s,a)\in \cS \times \cA} d_k(s,a)\cdot \nu(s,a) \le q_k, \ \forall\ k\in[K]. \notag
\end{align}

The second method utilizes the Lagrangian duality.
Let $\lambda \in \mathbb{R}^K$ denote the Lagrangian multipler. Define
\begin{align} \label{Lag}
L(\pi,\lambda):= C(\pi)+ \sum_{k=1}^K [\lambda]_k \cdot ( D_k(\pi)-q_k).
\end{align}
Then the CMDP can be equivalently formulated as 
$\inf_{\pi \in \Pi_S} \sup_{\lambda \ge 0} L(\pi,\lambda).$
By Theorem 3.6 in \cite{alt}, we can exchange the order of inf and sup and obtain,
%\begin{align}  \label{duality}
\[\inf_{\pi \in \Pi_S} \sup_{\lambda \ge 0} L(\pi,\lambda)=\sup_{\lambda \ge 0} \inf_{\pi \in \Pi_S}  L(\pi,\lambda) = \sup_{\lambda \ge 0} \inf_{\pi \in \Pi_D}  L(\pi,\lambda),\]
%\end{align}
where the last equation holds because for each fixed $\lambda$, the inner problem is an unconstrained MDP and the optimal policy is a stationary deterministic policy. 
We emphasize that given the optimal solution $\lambda^*$ to the dual problem, not every policy $\pi(\lambda^*)  $  that minimizes $ L(\pi,\lambda^*)$ is the optimal policy to the original CMDP. A necessary condition for $\pi(\lambda^*)$ to be optimal for the original CMDP is the complementary slackness: $ [\lambda^*]_k\cdot D_k(\pi(\lambda^*))=0, \forall k \in [K] $. %When the CMDP has only stochastic optimal policies, it is possible that there is no deterministic policy that minimizes $ L(\pi,\lambda^*)$ and satisfies the completeness and slackness condition simultaneously.  

The dual problem $\sup_{\lambda \ge 0} \inf_{\pi \in \Pi_D}  L(\pi,\lambda)$ leads to the following LP formulation:
\begin{align}  \label{dual_LP}
  &\max_{\phi,\lambda} \sum_{s\in \cS} \mu_0(s)\phi(s)-\sum_{k=1}^K[\lambda]_kq_k \\
  &\text{s.t.}\ \phi(s)\le (1-\gamma) \Big( c(s,a)+\sum_{k=1}^K[\lambda]_kd_k(s,a)\Big) + \gamma\cdot \sum_{s' \in \cS}\phi(s')P(s'|s,a), \notag
\end{align}
where $\phi(s)$ denotes the value function with initial state $s$. Note that \eqref{LP} and \eqref{dual_LP} are dual of each other.

Various methods have been developed in the literature to solve the LPs \eqref{LP} or \eqref{dual_LP}.
There are two main obstacles to solve the LPs in practice. First, it can be computationally prohibitive when dealing with a large state space or a large action space. Second, it requires explicit characterization of the transition kernel $P$. To overcome these difficulties, we next develop a sampling-based primal-dual algorithm to solve CMDPs. %Our algorithm builds on the Lagrangian defined in \eqref{Lag} .

\section{The Primal-Dual Algorithm}   \label{algorithm} 
Consider the Lagrangian dual problem
\begin{equation}\label{eq:dual}
 \sup_{\lambda \ge 0} \inf_{\pi \in \Pi_S} L(\pi,\lambda).
\end{equation}
For each fixed $\lambda$, the inner problem is an unconstrained MDP. A natural idea is to solve the unconstrained MDP via a sampling-based method and then update the Lagrangian multipliers via subgradient ascent. Such an idea is exploited in \citep{bpl}. However, this method is computationally expensive, since we need to solve a new MDP every time the Lagrangian multipliers are updated. %Multiple iterations are required to find the optimal policies with different Lagrangian multipliers. We would like to develop a provable algorithm that alternatively updates the policy and Lagrangian multiplier, instead of solving the penalized MDPs exactly. 
%Our method takes a similar primal-dual updating approach. The key difference is that instead of solving a new MDP at each iteration, we only do one regularized policy update at each iteration.
In contrast, our method only requires a single policy update at each iteration, i.e., we do not need to solve for the corresponding optimal policy at each iteration.

We develop the algorithm and analyze its convergence in  $\mathcal{M}(\Pi_S)$, the space of mixing stationary policies, rather than $\Pi_S$. The benefits of allowing the mixing are twofolds. First, it provides an intuitive way to understand strong duality:
\begin{align}  \label{strong_dual}
\inf_{\pi \in \cM({\Pi_S})}\sup_{\lambda \ge 0}L(\pi,\lambda)=\sup_{\lambda \ge 0}\inf_{\pi \in \cM({\Pi_S})}L(\pi,\lambda).
\end{align}
With the mixing operation, we can treat $C(\pi)$ and $D(\pi) $ as infinite-dimensional linear functions with respect to the distributions of initial randomization of  policies in $\Pi_S$. Hence, the Lagrangian $L(\pi,\lambda) $ is a bilinear function and strong duality follows from the minimax theorem \citep{Sion}. Second, in primal-dual algorithms, we in general need to take the average of the trajectories to obtain convergence \citep{ssp}. In our case, caution needs to be taken when defining the average. In particular, note that the objective and constraints are inner products of the cost functions and the occupation measures. Thus, what we need to average across are the occupation measures. However, since the mapping from the policy to the corresponding occupation measure is nonlinear, we cannot average the policy $\pi(\cdot|s)$, i.e., the probability of taking each action at each state, directly. The mixing operation provides a simple way to average the occupation measures. In addition, given a mixing policy, under mild regularity conditions,  there exists a non-mixing stationary policy that has the same occupation measure (Theorem 3.1 of \cite{alt}). In particular, for $\pi \in \cM(\Pi_S)$, let $\nu^{\pi}(\cdot,\cdot) $ be the corresponding occupation measure. Then, we can construct such a stationary policy $\tilde{\pi} $ via 
\begin{align} \label{trans}
  \tilde{\pi}(a|s)=\frac{\nu^{\pi}(s,a)}{\sum_{a\in \cA} \nu^{\pi}(s,a)}.
\end{align}

Our algorithmic development is based on strong duality \eqref{strong_dual},  which holds under certain regularity conditions (see Section \ref{theory} for details).
By the minimax theorem, there exists a saddle point $(\pi^*,\lambda^*)$ such that
\begin{align} \label{saddle_point}
L(\pi^*,\lambda) \le L(\pi^*,\lambda^*) \le L(\pi,\lambda^*),\ \forall\ \lambda \in \R^K_+, \pi\in \cM({\Pi_S}).
\end{align}
Moreover, $\pi^*$ is an optimal solution to the primal problem and $\lambda^* $ is an optimal solution to the dual problem. 
In addition, $L(\pi^*,\lambda^*)$ equals to the optimal cost of the CMDP.  The saddle point property \eqref{saddle_point} suggests that we can use iterative primal-dual updates to find the saddle point. %of \eqref{strong_dual}.

We next introduce our actual algorithm. Note that for a fixed value of $\lambda$, the inner inf-problem is an unconstrained MDP with modified instantaneous cost  $c^\lambda(s,a):=c(s,a)+\sum_{k=1}^K [\lambda]_k(d_k(s,a)-q_k)$. In what follows, we refer to the inner problem $\inf_{\pi \in \cM({\Pi_S})}L(\pi,\lambda)$ as the modified unconstrained MDP.

 For a given policy $\pi$ and Lagrangian multiplier $\lambda$, define
\begin{align} \label{Q_function}
Q^{\pi, \lambda} (s,a):=(1-\gamma)\cdot \Big( c^\lambda(s,a)+\E^\pi \Big[\sum_{t=1}^\infty \gamma^t c^\lambda(s_t,a_t)\big| s_0=s, a_0=a\Big]\Big),
\end{align}
which is known as the action-value function or $Q$-function. 
Let $\pi_m$ and $\lambda_m$ denote the policy and the Lagrangian multiplier obtained at iteration $m$. For the policy update, we use KL divergence as the regularization \citep{rmdp}. 
In particular, the regularized policy iteration is defined as
\begin{equation}\label{eq:iterate_policy}
\pi_{m}(a|s)=\argmin_{\pi(\cdot|s) \in \Delta_{\cA}}\Big\{ \big\langle Q^{\pi_{m-1},\lambda_{m-1}}(s,\cdot), \pi(\cdot|s) \big\rangle + \eta_{m-1}^{-1}\cdot \text{KL}\big(\pi(\cdot|s) \| \pi_{m-1}(\cdot|s)\big) \Big\}, 
\end{equation}
where $ \eta_{m-1}>0$ is the stepsize that determines the power of regularization. 
%Note that the regularized policy iteration takes the same form as the update of the mirror descent algorithm in its proximal form. 
Note that the regularized policy iteration  \eqref{eq:iterate_policy} is defined state-wise, i.e., for each $s\in \cS$. The minimization is taken over the probability simplex $\Delta_{\cA}:=\{  \pi(\cdot|s) : 0\le \pi(a|s)\le 1, \sum_{a \in \cA}  \pi(a|s) =1 \}$.

Let $\Lambda_M$ denote a suitably bounded domain that includes the dual optimal solution $\lambda^*$. We will provide an explicit construction of $\Lambda_M $  in Section \ref{theory}. 
To update the Lagrangian multiplier, we use the projected  subgradient ascent:
\begin{equation} \label{eq:iterate_lambda}
\lambda_{m} = \text{Proj}_{\Lambda_M} \Big\{ \lambda_{m-1}+ \eta_{m-1} \cdot \partial_\lambda L(\pi_{m-1},\lambda_{m-1}) \Big\},
\end{equation}
where $\text{Proj}_{\Lambda_M}\{\cdot\}$ denotes the projection (in $L^2$-norm) on $\Lambda_M$. We need such a projection to ensure the boundedness of ``subgradient" in order to establish convergence.  

By the definition of KL-divergence, the regularized policy iteration can be re-written as 
\begin{align} \label{policy_iteration}
\pi_{m}(\cdot|s) &= Z_{m-1}^{-1}\cdot \pi_{m-1}(\cdot|s)\cdot \exp \big\{-\eta_{m-1}\cdot Q^{\pi_{m-1},\lambda_{m-1}}(s,\cdot) \big\},
\end{align}
where $Z_{m-1}$ is some normalization constant. For the subgradient ascent update, we have
\begin{equation} \label{lambda_iteration}
\big[\partial_\lambda L(\pi_{m-1},\lambda_{m-1})\big]_k=D_k(\pi_{m-1})-q_k.
\end{equation}
Both \eqref{policy_iteration} and \eqref{lambda_iteration} can be evaluated/approximated using simulation.
More advanced approximation techniques for policy evaluation like TD-learning can also be applied here.

Suppose that our algorithm runs $T-1$ iterations and generates a sequence $\{(\pi_m,\lambda_m)\}_{0 \le m \le T-1}$. Then, the algorithm outputs a mixing policy and Lagrangian multiplier by taking a weighted average of the outputs:
\begin{align} \label{output}
 \bar{\pi}_T=\sum_{m=0}^{T-1} \tilde{\eta}_m  \pi_m,\ \bar{\lambda}_T=\sum_{m=0}^{T-1} \tilde{\eta}_m \lambda_m, \mbox{ where $\tilde{\eta}_m=\eta_m/\sum_{m=0}^{T-1} \eta_m$.}
\end{align}
The averaging operation is required for convergence, since the objective  $L(\pi,\lambda) $  is bilinear and does not possess sufficient convexity. In particular, counter-examples that fail to converge without averaging exist. The summation in the definition of $\bar{\pi}_T $ is interpreted as the mixing operation, i.e., it mixes policies $(\pi_{0},\ldots,\pi_{T-1})$ with initial randomization distribution $(\tilde{\eta}_0,\cdots,\tilde{\eta}_{T-1})$. 
Note that this essentially takes the average of the occupation measures of $\pi_m$'s. 
From $\bar\pi_T$, we can apply \eqref{trans} to define a non-mixing stationary policy that has the same occupation measure.  

Above all, our primal-dual algorithm is summarized in Algorithm \ref{main_algorithm}. 

\begin{algorithm}       
  \caption{Primal-Dual Algorithm to CMDP}
  \label{main_algorithm}
  \begin{algorithmic}
  \STATE {Input: pre-specified projection domain $\Lambda_M$, stepsizes $\{\eta_m\}_{m\geq 0}$, initial policy $\pi_0$ and Lagrangian multiplier $\lambda_0$} 
  \FOR{ $m=1,\ldots, T-1$ } 
  \STATE {update Lagrangian multipliers and policy as
    \begin{align*}
    \begin{cases}
    \lambda_{m} = \text{Proj}_{\Lambda_M} \Big\{ \lambda_{m-1}+ \eta_{m-1} \cdot \partial_\lambda L(\pi_{m-1},\lambda_{m-1}) \Big\}, \\
      \pi_{m}(\cdot|s) \propto \pi_{m-1}(\cdot|s)\cdot \exp \big\{-\eta_{m-1}\cdot Q^{\pi_{m-1},\lambda_{m-1}}(s,\cdot) \big\}.
     \end{cases}  
    \end{align*}
   }
  \ENDFOR
  \STATE{Output: mixing policy $\bar{\pi}_T=\sum_{m=0}^{T-1}\tilde{\eta}_m\pi_m $, where $\tilde{\eta}_m={\eta}_m/\sum_{m=0}^{T-1} {\eta}_m$.    }
  \end{algorithmic}  \label{main_alg}
\end{algorithm}

%Finally, note that by \eqref{policy_iteration}, when stepsize $\eta_m$ is large, the regularization effect of KL divergence is neglectable. Then the randomized action can also been approximated by greedy and hence, the regularized policy iteration reduces to the usual policy iteration. From this perspective, our algorithm executes policy iteration to improve policies and gradient ascent to find better Lagrangian multiplier, alternatively, instead of solving the unconstrained MDP exactly for each multiplier. As a result, the computational complexity is significantly improved compared with previous works. 

\section{Convergence Analysis} \label{theory}
In this section, we conduct detailed performance analysis of Algorithm \ref{main_algorithm}. 
%Recall that Algorithm \ref{main_algorithm} outputs a mixing policy 
%$$
%\bar{\pi}_T=\sum_{m=1}^T \tilde{\eta}_m \pi_m  \in \cM(\Pi_{S})
%$$ after $T$ iterations, where $\tilde{\eta}_m=\eta_m/\sum_{m=1}^T\eta_m $. 
In particular, we study the performance of policy $\bar{\pi}_T$ by analyzing the values of the objective $C(\bar{\pi}_T)$ and the constraints $D(\bar{\pi}_T)$. We show that the objective value $C(\bar{\pi}_T)$ converges to the optimal $C^*:=C(\pi^*)=L(\pi^*,\lambda^*)$ at a rate of $O(\log(T)/\sqrt{T})$. In addition, even though $\bar{\pi}_T $ may be infeasible, we show that the violation of constraints, which is measured by 
\begin{align} \label{violation}
\big\| [D(\bar{\pi}_T)-q]^+ \big\|:=\left(\sum_{k=1}^K \big([D_k(\bar{\pi}_T)-q_k]^+\big)^2\right)^{1/2},
\end{align}
converges to zero at a rate of $O(\log(T)/\sqrt{T})$. The analysis builds on a combination of subgradient method  for saddle point problem and mirror descent  for regularized policy iteration. %In the next, we first introduce some necessary assumptions in \ref{ass} and discuss their relevance. Then, we present the convergence results and proofs in Sections \ref{convergence} and \ref{proof}, respectively.

Recall that our algorithmic development builds on the strong duality of CMDP.  For CMDPs with finite state and action spaces, the strong duality always holds (Theorem 3.6 in \citep{alt}). However, when the state space is countably infinite, we need more regularity conditions to ensure the strong duality. One sufficient condition is that the instantaneous costs of the CMDP are uniformly bounded from below (see Definition 7.1, Theorem 9.9, and Chapter 10.3 in \citep{alt}). Specifically, we impose the following assumption. 
\begin{assumption} \label{ass0}
 [Lower Bound of Instantaneous Costs]
There exists a constant $W$ such that for all $s\in \cS$, $a \in \cA$, and $k=1,2,\ldots,K$, 
 $$
 c(s,a) > W, ~~d_k(s,a) >W. 
 $$
\end{assumption}

To establish the convergence result, we also require the Slater's condition: 
\begin{assumption} \label{ass1}
  [Slater's Condition] There exists some policy $\tilde{\pi}$ such that 
  $$
  D_k(\tilde{\pi}) < q_k,\ \forall 1\le k\le K.     
  $$
\end{assumption}
Slater's condition ensures the existences of finite and bounded optimal Lagrangian multipliers
$$
\lambda^*=\argmax_{\lambda \ge 0}\Big\{\inf_{\pi \in \cM(\Pi_S)} L(\pi,\lambda)   \Big\}.
$$
This condition is commonly assumed in the constrained optimization literature. For many practical problems, the Slater's condition holds trivially. 

Our last assumption is about the boundedness of the ``subgradient'', which regularizes the movement of policies and Lagrangian multipliers at each iteration. 
Recall that in Algorithm \ref{main_algorithm}, after applying the subgradient ascent for Lagrangian multipliers, we project  $\lambda $ onto a bounded domain $\Lambda_M $, which  takes the form
\begin{align} \label{Lambda_area}
\Lambda_M=\big\{ \lambda \in \R^K_+: \| \lambda \| \le M+r \big \},
\end{align}
where $M$ is an upper bound of $\| \lambda^*\| $ and $r>0$ is a slackness constant. %We impose the following assumption.  

\begin{assumption} \label{ass2} [Bounded Subgradient] There exists some constant $G>0$ such that for any $\lambda \in \Lambda_M$ and policy $\pi \in \cM(\Pi_S),$  
$$\big\|\partial_{\lambda} L(\pi,\lambda)  \big\| \le G,\ \sup_{s\in \cS}\sup_{a\in \cA}  \big|Q^{\pi,\lambda}(s,a)\big| \le G. $$
\end{assumption}

Since $Q^{\pi,\lambda}(s,a)$ is linear in $ \lambda$, it is necessary to restrict $\lambda$ to a bounded domain $\Lambda_M$ for Assumption \ref{ass2} to hold. That is  why we need the projection step in updating $\lambda$. Note that when the instantaneous cost functions $c(\cdot,\cdot) $ and  $d_k(\cdot,\cdot)$ are uniformly bounded or when the state and action spaces are finite, Assumption \ref{ass2} holds trivially.   

Lastly, we comment that the Slater's condition (Assumption \ref{ass1}) not only guarantees the existence and boundedness of  $\lambda^*$, but also provides an explicit upper bound of $\|\lambda^*\| $. In particular, let $\tilde{\pi}$ be a Slater point (a policy that satisfies the Slater's condition), then we have  
\begin{align} \label{LM_upper}
\|\lambda^*\| \le  -\frac{C(\tilde{\pi})- \tilde{c}}{\max_{1\le k \le K}\big\{D_k(\tilde{\pi})-q_k\big\}},
\end{align}
where $\tilde{c} \le C(\pi^*)$ is an arbitrary lower bound for the dual problem.  In many applications, it is possible to obtain a better upper bound of $\| \lambda^*\| $ than \eqref{LM_upper} by exploiting the structure of the specific problem. 
%Note that the upper bound $G$ in Assumption \ref{ass2} depends on the shape of $ \Lambda_M$ obviously. As we will show later, $ G$ influences the convergence rate of our algorithm  the upper estimate of optimal Lagrangian multiplier $M$ enters the convergence rate through $G$. Hence, a more accurate estimate of $\|\lambda^*\| $ means more efficiency.

%\subsection{Convergence Result} \label{convergence}
%We present the convergence result of our algorithm in this section. We show that under Assumptions \ref{ass1} and \ref{ass2}, with suitable choices of stepsize, the cost of  policy $ \bar{\pi}_T$ converges to the optimal at a rate of $O(\log(T)/\sqrt{T})$ or to the neighborhood of the optimal at a rate of $O(1/T)$. Similarly, the violation of constraints, which is is measured via \eqref{violation}, converges to zero or to the neighborhood of zero, at the same speed. 

Next, to establish the convergence, we need to construct an appropriate potential function, which is also known as Bregman divergence in the optimization literature.
The potential function ensures that the regularized policy iteration is equivalent to minimize the sum of a linear approximation of the objective function and the potential function. We next introduce this potential function, which is essentially a weighted KL-divergence.
 
Consider the state occupation measure $\nu_s^{\pi}$ induced by a policy $\pi \in \Pi_S$, i.e., 
$$
\nu_s^\pi(s):=(1-\gamma)\cdot \E_{s_0 \sim \mu_0} \Big[ \sum_{t=0}^\infty \gamma^t P^\pi(s_t=s | s_0) \Big].
$$
The KL-divergence between two stationary policies $\pi_1$ and $\pi_2$ weighted by $\nu_s^{\pi} $ is defined as     
\begin{align} \label{w-KL}
  \Phi^{\pi}(\pi_1\| \pi_2)=\E_{s\sim \nu_s^{\pi}}\Big[\KL\big(\pi_1(\cdot|s)\| \pi_2(\cdot|s)\big)\Big].
\end{align}
When $\pi_1$ and $\pi_2$ are mixing policies, we first transform them to the equivalent stationary policies via \eqref{trans}, and then define $ \Phi^{\pi}(\pi_1\| \pi_2)$ as the weighted KL-divergence between the equivalent stationary policies.

By definition, $\Phi^{\pi}(\pi_1\| \pi_2)$ measures the discrepancy between two policies weighted by a given state occupation measure. It connects the regularized policy iteration in \eqref{eq:iterate_policy}, which is defined state-wise, with a single objective and serves as the Bregman divergence in mirror descent analysis. Unlike the traditional analysis of mirror descent where the potential function is fixed \citep{md}, in the analysis of regularized policy iteration, we need to construct a policy-dependent potential function
and cannot fix the weight of KL-divergence. 
However, since policy updates are defined state-wise, for an arbitrary weight, the regularized policy iteration always takes the form of minimizing a linear approximation of the objective function regularized by a certain potential function. Thus, the analysis of mirror descent can be applied here with some modifications.  

We are now ready to introduce the convergence results of our primal-dual algorithm. 
\begin{theorem} \label{thm1} (Convergence of Main Algorithm)
  Under Assumptions \ref{ass0}-\ref{ass2}, if the step size $\eta_m=\Theta(1/\sqrt{m})$, then there exist positive constants $\kappa_1$ and $\kappa_2$, such that 
  \begin{align*}
  \big\| [D(\bar{\pi}_T)-q]^+ \big\|  \le  \Big(G^2   \big( 1+ \frac{5}{8}\kappa_2\log(T) \big)+ \Phi^{\pi^*}(\pi^*\| \pi_0) \Big )\frac{1}{2r(1-\gamma)\kappa_1  \sqrt{T}},
  \end{align*}
and 
\begin{align*} 
 C(\bar{\pi}_T)-L(\pi^*,\lambda^*) & \le  \Big( \frac{5G^2}{8} \cdot \kappa_2 \cdot \log(T) +\Phi^{\pi^*}(\pi^*\| \pi_0) +\frac{\| \lambda_0 \|^2}{2} \Big)\frac{1}{(1-\gamma)\kappa_1 \sqrt{T}},\\
 C(\bar{\pi}_T)-L(\pi^*,\lambda^*)  &\ge - \| \lambda^* \| \Big(G^2   \big( 1+ \frac{5}{8}\kappa_2 \log(T)\big)+ \Phi^{\pi^*}(\pi^*\| \pi_0) \Big )\frac{1}{2r(1-\gamma)\kappa_1 \sqrt{T}}.
\end{align*}
If the step size is constant $\eta_m=\eta$, then 
\begin{align*}
 \big\| [D(\bar{\pi}_T)-q]^+ \big\| \le  \big(G^2 +  (1-\gamma)^{-1} \cdot \Phi^{\pi^*}(\pi^*\| \pi_0) \big ) \frac{1}{2r T\eta} +  \Big(\frac{1}{2}+\frac{1}{8(1-\gamma)}\Big)\frac{G^2\eta}{2r},
\end{align*} 
and
\begin{align*}
C(\bar{\pi}_T)-L(\pi^*,\lambda^*)  \le&   \left(  (1-\gamma)^{-1}\cdot \Phi^{\pi^*}(\pi^*\| \pi_0) +\frac{\| \lambda_0 \|^2}{2} \right)\frac{1}{T\eta} + \frac{5G^2\eta}{8(1-\gamma)}, \\
C(\bar{\pi}_T)-L(\pi^*,\lambda^*) \ge& -\| \lambda^*\|  \big(G^2 +  (1-\gamma)^{-1} \cdot \Phi^{\pi^*}(\pi^*\| \pi_0) \big ) \frac{1}{2r T\eta}
-\| \lambda^*\|\Big(\frac{1}{2}+\frac{1}{8(1-\gamma)}\Big)\frac{G^2\eta}{2 r},
\end{align*}
where $r$ is the slackness constant  in \eqref{Lambda_area}.
\end{theorem}

Theorem \ref{thm1} indicates that with decreasing step size, $\eta_m=\Theta(1/\sqrt{m})$, our primal-dual algorithm achieves $O(\log(T)/\sqrt{T})$ convergence. In particular, 
\[
 \big\| [D(\bar{\pi}_T)-q]^+ \big\| = O(\log(T)/\sqrt{T}) \mbox{ and } |C(\bar{\pi}_T)-L(\pi^*,\lambda^*)|=O(\log(T)/\sqrt{T}).
\]
For constant step size, $\eta_m=\eta$, our primal-dual algorithm converges to a neighborhood of the optimal at rate $O(1/T)$. In particular,
\[
 \big\| [D(\bar{\pi}_T)-q]^+ \big\| = O(1/(\eta T)+\eta) \mbox{ and } |C(\bar{\pi}_T)-L(\pi^*,\lambda^*)|=O(1/(\eta T)+\eta)
\]
These convergence rates match those in \cite{bpl}, which requires solving the modified unconstrained MDP to the optimal at each iteration.
We also note that it is unlikely to improve the convergence rate beyond {$\Theta(1/\sqrt{T})$}. This is because the dual problem is a finite-dimensional concave optimization problem without strong concavity. The convergence rate of the subgradient method in this case is lower bounded by $\Omega(1/\sqrt{T})$ \citep{bubeck}. The proof of Theorem \ref{thm1} is deferred to the appendix. 

We comment that in the bounds in Theorem \ref{thm1}, although the slackness constant $r$ appears in denominators only, the constant $G$, which is an upper bound of the subgradients,  grows linearly in $r$. In particular, by Assumption \ref{ass2}, $G$ is determined by the shape of $ \Lambda_M$.  Hence, $r$ cannot be set too large.
%Hence, there is a tradeoff in the choice of $r$. %See \citep{ssp} for a discussion on the optimal choice of $r$. }

%In Theorem \ref{thm1}, if we ignore all the constants, then the violation of constraints and accumulated cost of policy $\bar{\pi}_T$ converge to zero and the optimal, respectively, at a rate of $O(\log(T)/\sqrt{T})$ with the $O(1/\sqrt{m})$ stepsize.  They also converge to the neighborhood of zero and the optimal at a rate of $O(1/T)$ with constant stepsize. Such a rate of convergence matches the results in \citep{bpl}, where they need multiple iterations to solve the penalized MDPs with every Lagrangian multiplier. However, our primal-dual algorithm only executes a single (regularized) policy iteration for each Lagrangian multiplier, which reduce the computational complexity significantly. If we ignore the $O(\log(T))$ factor, the $O(1/\sqrt{T})$ are unlikely to be improved, since the dual problem is just a finite dimensional concave optimization problem without strongly concavity guarantee. It is well-known that the lower bound of the convergence rate of the subgradient method in this case is $O(1/\sqrt{T})$. Meanwhile, since regularized policy iteration has the similar structure of mirror descent update, only the $O(1/\sqrt{T})$ rate of convergence is available for unconstrained MDPs. See, for example, \citep{Liu} for more regards on these analysis. 

\section{Weakly Coupled MDP and Weakly Coupled CMDP}  \label{WCMDP}
One fundamental challenge in solving MDPs and CMDPs is the curse of dimensionality. 
%{\color{blue} As the dimension of the problem increases, the complexity of most algorithms increases exponentially with the dimension. ?}
However, there is an important class of problems that has certain decomposable structures.
These problems, which are often referred to as weakly coupled MDPs/CMDPs, contain multiple subproblems which are almost independent of each other except for some linking constraints on the action space \citep{Singh}. More precisely, 
for a weakly coupled MDP consisting of $I$ sub-problems $ \{(\cS^i,\cA^i,P^i,c^i(\cdot,\cdot),\gamma,\mu_0^i)\}_{i \in [I]}$, we have the following structural properties:
P1. Its state and action spaces can be expressed in the form of Cartesian products, i.e.,
\begin{align*}
  \bm{s}&=(s^1,\ldots,s^I),\ \bm{\cS}=\cS^1 \times \cS^2 \times \ldots \times \cS^I, \\
  \bm{a}&=(a^1,\ldots,a^I),\ \bm{\cA}=\cA^1 \times \cA^2 \times \ldots \times\cA^I.
\end{align*}    
P2. For each state $\bm{s}_t $ and action $\bm{a}_t$, the instantaneous cost admits an additively separable form  
$$c(\bm{s}_t,\bm{a}_t)=\sum_{i=1}^I c^i(s_t^i,a_t^i).$$ 
P3. The joint initial distribution satisfies $\bm{\mu}_0(\bm{s})=\mu^1_0(s^1) \cdot \mu^2_0(s^2) \cdot \ldots \cdot \mu^I_0(s^I)$ and the one-step transition dynamics of the sub-MDPs are independent of each other, i.e., 
\begin{align*} 
P(\bm{s}_{t+1}|\bm{s}_{t},\bm{a}_t)=\prod_{i=1}^I P^i(s^i_{t+1}|s_t^i,a_t^i). 
\end{align*}
For the linking constraints, let  $b^i(\cdot,\cdot): \cS^i\times\cA^i \to \mathbb{R}^K$ be a $K$-dimensional real function, which can be interpreted as the resource consumed by the $i$-th sub-problem, $i\in [I]$. 
Then, at each state  $\bm{s}$, the feasible actions need to satisfy
\begin{align}  \label{joint_constraint}
 % \bar{\cA}(\bm{s})=\Big\{\bm{a}=(a^1,\ldots,a^I) \in \bm{\cA}:\  \sum_{i=1}^I b^i(s^i,a^i) \le \textbf{0} \Big\}. 
b(\bm{s},\bm{a})= \sum_{i=1}^I b^i(s^i,a^i) \le q.
\end{align} 
where $q \in \mathbb{R}^K$. Note that the linking constraint \eqref{joint_constraint} is a hard constraint and needs to be satisfied path-by-path almost surely. 
For a weakly couple CMDP, it satisfies the same structural properties, P1-P3, as the weekly coupled MDP.
The only difference is that the linking constraint now takes the form
\begin{align} \label{joint_constraint2}  
(1-\gamma)\cdot  \E_{\bm{s}_0\sim \bm{\mu}_0}\Big[  \sum_{t=0}^\infty   \gamma^t \cdot \sum_{i=1}^I  b^i(s^i_t,a^i_t)  \big| \bm{s}_0 \Big] \le q.
 \end{align} 
 
The weakly coupled MDP and the weakly coupled CMDP are closely related to each other. 
%In particular, in CMDP, we relax the ``hard" path-by-path linking constraint \eqref{joint_constraint} in the corresponding MDP to a ``soft" constrain on the expected discounted average \eqref{joint_constraint2}. We next discuss the connection between the two in more details.
%Compared with the weakly coupled MDP, the weakly coupled CMDP first incorporates the path-wise constraints into a single objective and then replaces the almost sure constraints to expectation constraints. So the weakly coupled MDP is essentially an unconstrained MDP (although the hard constraints determine the shape of action spaces) and weakly coupled CMDP is a standard CMDP (with expectation type constraints). 
Let
$$\bar{\cA}(\bm{s})=\Big\{\bm{a}=(a^1,\ldots,a^I) \in \bm{\cA}:\  \sum_{i=1}^I b^i(s^i,a^i) \le q \Big\} $$
be the (joint) action space of  a weakly coupled MDP. Then, the Bellman equation  is 
\begin{align*}
V^{*}_{\bm{\mu}_0}(\bm{s})=\min_{ \bm{a} \in \bar{\cA}(\bm{s})}\Big\{\sum_{i=1}^I c^i(s^i,a^i) +\gamma\cdot \sum_{\bm{s}'\in \bm{\cS}} V^{*}_{\bm{\mu}_0}(\bm{s}') \cdot P(\bm{s}'|\bm{s},\bm{a})  \Big\}.
\end{align*}
When the number of sub-MDPs $I$ is large, even if the scale of each subproblem is small, the size of joint state space $\bm{\cS}$ can be prohibitively large. Hence, solving the MDP directly can be intractable. 
Two decomposition schemes have been proposed to alleviate the curse of dimensionality: LP-based ADP relaxation and Lagrangian relaxation \citep{Adelman}. Both of them lead to $I$ independent sub-LPs, which reduces the complexity significantly. 
The LP-based ADP relaxation approximates the value function with additively separable functions,  
i.e., $V^{*}_{\bm{\mu}_0}(\bm{s}) \approx \sum_{i=1}^I V^*_{{\mu}^i_0}(s^i).$
The Lagrangian relaxation dualizes the constraints \eqref{joint_constraint} based on the LP representation of the Bellman equation.
The latter relaxation translates the weakly coupled MDP to a weakly coupled CMDP. It has been established that the 
optimal cost of the relaxed CMDP provides a lower bound for the optimal cost of the original MDP \citep{Adelman}. 
%We also comment that the policy derived based on the CMDP may violate the hard constraints in the corresponding MDP. 
%Several techniques have been proposed to reduce the probability of constraint violation, see \citep{Balseiro} for an example. 

Many Operations Management problems can be formulated as weakly coupled MDPs/CMDPs. Examples include inventory planning problems with multiple types of inventories and budget constraints, and scheduling of parallel-server queues with multiple classes of customers. We provide more details about these problems in Sections \ref{sec:inventory} and \ref{sec:queue}, where we apply our primal-dual algorithms to solve them.

When applying the primal-dual algorithm to solve weakly coupled CMDPs, it can be easily adapted to enjoy the decomposability. 
We call a policy $\bm{\pi}$ decomposable if it takes the product form: 
\[\bm{\pi}(\bm{a}|\bm{s})=\prod_{i=1}^I \pi^i(a^i|s^i).\]
Since our algorithm converges with any initial policy, we shall start with a decomposable policy.
%which makes further update on policies also decomposable and simplifies the implementation of algorithm.
Let $\{\bm{s}_t\}_{t\ge 0}=\{(s^1_t,\ldots,s^I_t)\}_{t\ge 0}$ and  $\{\bm{a}_t\}_{t\ge 0}=\{(a^1_t,\ldots,a^I_t)\}_{t\ge 0}$ be the trajectory of the CMDP under policy $\bm{\pi}=(\pi^1,\ldots,\pi^I)$. To simplify the notations, for each $i\in [I]$, we define 
\begin{align*}  
  C^i(\pi^i)=(1-\gamma)\cdot\E^{\pi_i}_{s^i_0\sim \mu^i_0}\Big[  \sum_{t=0}^\infty \gamma^t\cdot c^i(s^i_t,a^i_t) \big| s^i_0 \Big],\
  B^i(\pi^i)=(1-\gamma)\cdot\E^{\pi_i}_{s^i_0\sim \mu^i_0}\Big[  \sum_{t=0}^\infty \gamma^t\cdot b^i(s^i_t,a^i_t) \big| s^i_0 \Big].
\end{align*}
Then, the CMDP can be written as 
\begin{align} \label{relaxed_cmdp}
\min_{(\pi^1,\ldots,\pi^I)} \sum_{i=1}^I C^i(\pi^i),\quad \text{s.t.}  \sum_{i=1}^I B^i(\pi^i) \le q.
\end{align}
%Now we use the primal-dual algorithm to solve the CMDP \eqref{relaxed_cmdp}, starting from a decomposable initial policy. 
When applying the primal-dual algorithm, if we start with a decomposable policy, then the policies obtained in all subsequent iterations are decomposable.
To see this, we note that the Lagrangian function,
\begin{align*}
L(\bm{\pi},\lambda)=\sum_{i=1}^I \big( C^i(\pi^i)+\lambda^\top B^i(\pi^i) \big)    - \lambda^\top q ,
\end{align*}
can be decomposed into $I$ independent subproblems.
%By Algorithm \ref{algorithm}, we execute the iteration,
%\begin{align*}  
%\begin{cases}
%\lambda_{m+1} =\text{Proj}_{\Lambda_M} \big\{ \lambda_{m+1}+\eta_m\cdot \partial_\lambda L_{\bm{\mu}_0}(\bm{\pi},\lambda) \big\}\\
%\bm{\pi}_{m+1}(\cdot|\bm{s}) \propto \bm{\pi}_{m}(\cdot|\bm{s})\cdot \exp \big\{-\eta_m \cdot Q^{\bm{\pi}_m,\lambda}(\bm{s},\cdot) \big\},
%\end{cases}
%\end{align*} 
%where $Q^{\bm{\pi}_m,\lambda}$ is the $Q$-function for the joint CMDP with penalized instantaneous cost $\sum_{i=1}^I c^i(\cdot,\cdot) +\lambda^\top d^i(\cdot,\cdot)$. Since $\bm{\pi}_m$ is decomposable, by the independency of sub-MDPs, we have 
If $\bm{\pi}_m$ is decomposable, 
\begin{align*}
  Q^{\bm{\pi}_m,\lambda}(\bm{s},\cdot)=\sum_{i=1}^I  Q^{{\pi}^i_m,\lambda}(s^i,\cdot),
\end{align*}
where $Q^{{\pi}^i_m,\lambda}(\cdot,\cdot)$ is the $Q$-function of the $i$-th modified sub-MDP with  instantaneous cost $c^i(\cdot,\cdot) +\lambda^\top b^i(\cdot,\cdot)$. Here, we ignore the constant $\lambda^Tq$, since subtracting a common constant in the $Q$-function does not change the updates of regularized policy iteration. This indicates that the regularized policy iteration,  including policy evaluation and  improvement, can be implemented separately in parallel via 
$$
{\pi}^i_{m+1}(\cdot|s^i) \propto {\pi}^i_{m}(\cdot|s^i)\cdot \exp \Big\{-\eta_m \cdot Q^{{\pi}^i_m,\lambda}(s^i,\cdot) \Big\},\ \forall i \in [I].
$$
Moreover, as the subgradient of Lagrangian multiplier takes  form 
$\partial_{\lambda} L(\bm{\pi}_m,\lambda)=\sum_{i=1}^I B^i(\pi_m^i)-q$,
it can be evaluated for the sub-MDPs in parallel as well.
Above all, in this case, the primal-dual algorithm improves the computational complexity from depending exponentially on $I$ to linearly on $I$.

 \section{Application to an Inventory Planning Problem}   \label{sec:inventory}
In this section, we apply the primal-dual algorithm to solve a multi-product multi-period newsvendor problem with budget constraints.

Consider the inventory planning problem with $I$ distinct products. At the beginning of each period, we need to decide the quantities to order based on the current inventory levels. The orders are assumed to be fulfilled without delay.
After the inventory is replenished, a random demand is realized. We assume the demands for each product are independent. Let $F_i$ denote the cumulative distribution function of the demand for product $i$ in each period. In particular, for each period, the demand for product $i$ is an independent draw from the distribution $F_i$.
For each product $i \in [I]$, we denote its inventory level at the beginning of period $t$ by $s^i_t$, the quantity we ordered by $a^i_t$, and the demand in period $t$ by $w_t^i$. For product $i$ in period $t$, if the demand does not exceed the current inventory level, i.e.,  $w^i_t \le s^i_t+a^i_t$, all the demand is fulfilled and the remaining inventory can be carried to the next period. Otherwise, only $s^i_t+a^i_t$ units are fulfilled in the current period. The remaining $(w^i_t-s^i_t-a^i_t)$ units are carried to the next period as backlog. We allow $s^i_t$'s to be negative to represent backlogs. For product $i$, inventory incurs a holding cost of $h_i$ per unit per period and backlog incurs a backlog cost of $b_i$ per unit per period.
%We are also subject to some budget constraints. Suppose that 
In addition, product $i$ in inventory consumes $v_i$ resource per unit per period. For a fixed $q>0$, we impose the following budget constraint  
\begin{align} \label{inventory_constraint}
(1-\gamma)\cdot \E\Big[ \sum_{t=0}^\infty  \sum_{i =1}^I  \gamma^t \cdot [s^i_t+a^i_t]^+ \cdot v_i  \Big | (s^1_0,\ldots,s^I_0) \Big] \le q. 
\end{align}  
The resource can be interpreted as, for example, the volume of each product. In this case, the above constraint put restrictions on the warehouse space.

The inventory planning problem can be formulated as a weakly coupled CMDP with state $\bm{s}=(s^1,\ldots,s^I)$, action $\bm{a}=(a^1,\ldots,a^I)$, and transition dynamics 
\[s^i_{t+1}=s^i_{t}+a^i_{t}-w^i_t,\ w^i_t \sim F^i(w),\ \forall i\in[I].\]
As the demands are independent, $ P(\bm{s}_{t+1}|\bm{s}_{t},\bm{a}_{t+1} )=\prod_{i=1}^I P({s}^i_{t+1}|{s}^i_{t},{a}^i_{t+1} )$.
The instantaneous cost function and auxiliary cost function are 
\begin{align*}
c(\bm{s}_t,\bm{a}_t) &= \sum_{i=1}^I h_i\cdot [s^i_{t}+a^i_{t}-w^i_t]^+ + b_i\cdot [w^i_t-s^i_{t}-a^i_{t}]^+, \\
b(\bm{s}_t,\bm{a}_t)&=\sum_{i =1}^I [s^i_t+a^i_t]^+ \cdot v_i.
\end{align*}  
%Due to the weak coupling structure, we can solve the multi-product news-vendor problem with constraint using our primal-dual algorithm efficiently, even if the number of products is large. 

To verify the correctness of our convergence analysis, we consider a small-scale instance of the problem with appropriate truncations. Such a truncation makes the state and action spaces finite. 
In this case, the optimal cost can be solved numerically (using the LP formulation). 
In particular, consider $I=2$, and demands for the two products are both uniformly distributed on set $\{1,2,\ldots,10\}$,
We impose an upper bound $10$ and a lower bound $-10$ for the state space. In particular, 
when backlogs drop below $-10$, the excess demands are lost without incurring any cost.
For other systems parameters, we set the holding costs $h_1=1,h_2=2$, backlog costs $b_1=2,b_2=3$, resource consumptions $v_1=1.5,v_2=1 $, threshold $q=10$, 
and discount rate $\gamma=0.75$.  
%With these parameters, the numbers of states and actions are roughly $400$.

%The state and action spaces of two products are the same. To simplify implementation, we also use truncation to make the state space finite. We set an upper bound $U=10$ and a lower bound $L=-10$ for the state space. Specifically, when the demand satisfies $s^i_t+a^i_t-w^t<-10$, we just set next state $s^i_{t+1}=-10$. The truncated demand will be lost forever without incurring any cost. In this case, the backlogged cost of product $i$ is $ b^i\cdot (-10+s^i_t-a^i_t)$. The upper bound $U=10$ requires that $0 \le a^i_t \le 10-s^i_t$. We assume that all the distributions of demands are the uniform distribution on set $\{1,2,\ldots,10\}$. The discount rate is $\gamma=0.75$. With these parameters, the numbers of states and actions are roughly $400$. 

\begin{figure}[htbp]
\centering

\subfigure[$C(\pi_T)$]{
\begin{minipage}[t]{0.5\linewidth}
\centering
\includegraphics[width=3in]{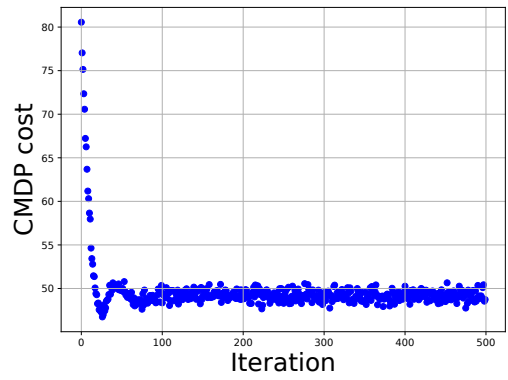}
%\caption{fig1}
\end{minipage}%
}%
\subfigure[$\sum_{t=0}^{T-1} \tilde{\eta}_t C(\pi_t)$]{
\begin{minipage}[t]{0.5\linewidth}
\centering
\includegraphics[width=3in]{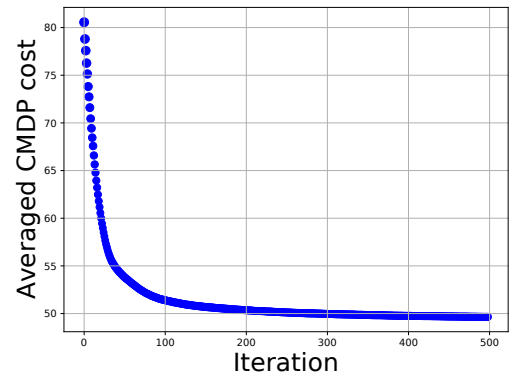}
%\caption{fig2}
\end{minipage}%
}%
                
\subfigure[{$\|[\sum_{i=1}^{2}B^i(\pi_T)-q]^+\|$}]{
\begin{minipage}[t]{0.5\linewidth}
\centering
\includegraphics[width=3in]{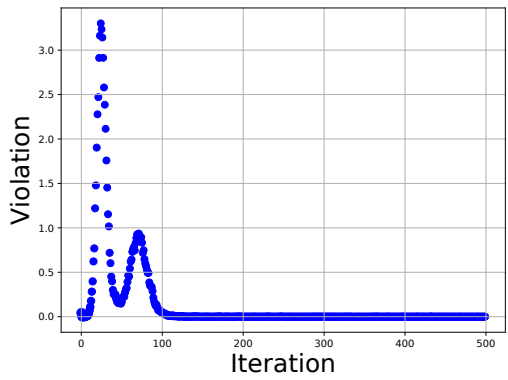}
%\caption{fig2}
\end{minipage}
}%
\subfigure[{$\sum_{t=0}^{T-1} \tilde{\eta}_t \|[\sum_{i=1}^{2}B(\pi_t)-q]^+\|$}]{
\begin{minipage}[t]{0.5\linewidth}
\centering
\includegraphics[width=3in]{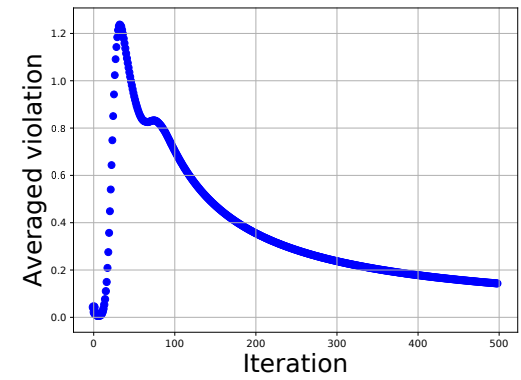}
%\caption{fig2}
\end{minipage}
}%
\centering
\caption{Trajectories of costs and constraints with constant step sizes }   \label{inv_fig1}
\end{figure}

\begin{figure}[htbp]
\centering
\subfigure[$\eta_m=0.2$ ]{
\begin{minipage}[t]{0.5\linewidth}
\centering
\includegraphics[width=3in]{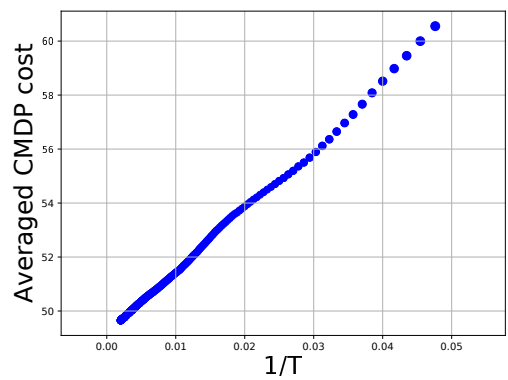}
%\caption{fig1}
\end{minipage}%
}%
\subfigure[$\eta_m=0.2/\sqrt{m+1}$]{
\begin{minipage}[t]{0.5\linewidth}
\centering
\includegraphics[width=3in]{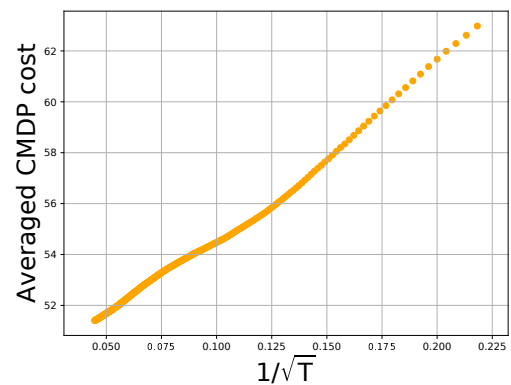}
%\caption{fig2}
\end{minipage}%
}%
\centering
\caption{Convergence rate with constant and decreasing step sizes}  \label{inv_fig2}
\end{figure}

When implementing the primal-dual algorithm, we use the standard Monte Carlo method to estimate the $Q$-function for a given policy. Since the system scale is small, we can enumerate all the state-action pairs in policy evaluation, i.e., no approximation of the value function is needed. 
The estimation of $Q$-function is based on an average of $400$ independent replications of the inventory process over $40$ periods of time.
%As $0.75^{40}\approx 10^{-5}$, the truncation is likely to have almost negligible bias. 
%Since the value function is the expectation of an infinite sum, when estimating $Q$, we truncate at the $40$-th period. As $0.75^{40}\approx 10^{-5}$, the truncation is likely to have almost negligible bias. 
%In each Monte Carlo estimation, we replicate $400$ times and take average. 
We implement two versions of the algorithm, one with constant step sizes $\eta_m=0.2$, the other with decreasing step size $\eta_m=0.2/\sqrt{m+1}$.
In each experiment, we run $500$ iterations in total and calculate the objective values for each iteration.    
%We compare the performance of our algorithm with the optimal cost calculated by solving the LP \eqref{dual_LP}, which is equal to $46.47$ based on our system parameters.  
%We check our theory in terms of objective value and the violation of constraint. 
%Let the policy obtained in the $m$-th iteration be $\pi_m $ and corresponding cost be $C(\pi_m)$.  The averaged CMDP cost by the $m$-th iteration is defined as $ \sum_{t=1}^m \tilde{\eta}_m C(\pi_m) $, which is same as the cost of the mixing policy $C(\bar{\pi}_m)=C(\sum_{t=1}^m \tilde{\eta}_m \pi_t)$ in expectation. The violation of constraint of policy $\pi_m$ is   
%$[B(\pi_m)]^+=[\E[\sum_{t=0}^\infty \gamma^t \sum_{i=1}^2 (v_i\cdot (s^i_t+a^i_t)-V ) ] ]^+$  and the averaged violation is $\sum_{t=1}^m \tilde{\eta}_m [B(\pi_m)]^+ $, which is an upper bound of the violation of mixing policy $\bar{\pi}_m$. 
The results of numerical experiments are summarized in Figures \ref{inv_fig1} and \ref{inv_fig2}. Figure \ref{inv_fig1} shows the trajectories of objective values and the constraint violations for different iterations with constant step size. 
We observe that after $500$ iterations, the averaged CMDP cost (without multiplying the $(1-\gamma)$ factor) converge to $49.26$, which is close to the optimal value $46.47$. In terms of feasibility, we calculate the violation of constraints, which is the expected value of the auxiliary cost minus the budget threshold. We observe that the averaged violation value converges to $0.1$ and many policies in the last iterations do not violate the constraint at all. Figure  \ref{inv_fig2} shows the relationship between $ \sum_{t=0}^{T-1} \tilde{\eta}_t C(\pi_t) $ and the reciprocal of the number of iterations  (for constant step size) or the reciprocal of square root of the number of iterations (for decreasing step size). In both cases, we observe a straight line, which confirms the rates of convergence developed in Theorem \ref{thm1}. 
%{\color{blue} Yi: what is Violation? Is the the value of the auxiliary cost or the probability of violating the constraints?}

\section{Application to Queueing Scheduling} \label{sec:queue}

In this section, we apply our primal-dual algorithm to a queue scheduling problem, which is motivated by applications in service operations management. Service systems often feature multiple classes of customers with different service needs and multiple pools of servers with different skillsets. Efficiently matching customers with compatible servers is critical to the management of these systems. In this context, we consider a parallel-server system (PSS) with multiple classes of customers and multiple pools (types) of servers. Customers waiting in queue incur some holding costs and routing customers to different pools leads to different routing costs. The goal is to find a scheduling policy that minimizes the performance cost (holding cost plus routing cost). This class of problems is known as the skill-based routing problem and has been widely studied in the literature. We refer to \citep{queue_survey} for a comprehensive survey of related works. 

In what follows, we first introduce the queueing model and some heuristic policies adapted from policies developed in the literature. We then present the implementation details of our primal-dual algorithm in this setting. Due to the large state and action spaces, we combine our primal-dual algorithm with several approximation techniques. Lastly, 
we compare the performance of our policy with the benchmark policies numerically. 

\subsection{Model and Benchmarks} 
\begin{figure}[htbp]
\centering
\includegraphics[width=2.75in]{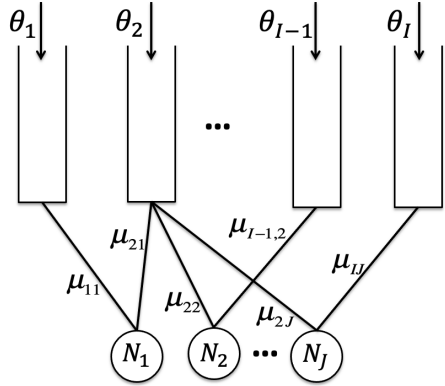}
\caption{Multi-class multi-pool queueing system}   \label{dynamics}
\end{figure}

The multi-class multi-pool queuing network has $I$ classes of customers and $J$ pools of servers. 
We consider a discrete time model. In each period, the number of arrivals of class $i$ customers follows a Poisson distribution with rate $\theta_i$.
There are $N_j$ homogeneous servers in pool $j$, $j \in [J]$. We assume that each customer can only be served by one server and each server can only serve one customer at a time. If a class $i$ customer is served by a server from pool $j$, its service time follows a geometric distribution with success probability $\mu_{ij}$. 
When there is no compatibility between customer class $i$ and server type $j$, $\mu_{ij}=0$. Figure \ref{dynamics} provides a pictorial illustration of such a system. 

%\begin{figure}[htbp]
%\centering
%\subfigure[multi-class multi-pool service system ]{
%\begin{minipage}[t]{0.5\linewidth}
%\centering
%\includegraphics[width=2.75in]{dynamics.eps}
%%\caption{fig1}
%\end{minipage}%
%}%
%\subfigure[decoupled subproblem]{
%\begin{minipage}[t]{0.5\linewidth}
%\centering
%\includegraphics[width=3.2in]{single_dynamics.eps}
%%\caption{fig2}
%\end{minipage}%
%}%
%\centering
%\caption{Dynamics of queueing model}   \label{dynamics}
%\end{figure}

We consider non-preemptive scheduling policies.
Let $A_i(t)$ denote the number of new class $i$ arrivals in time period $t$, i.e., $A_i(t)$ follows a Poisson distribution with rate $\theta_i$. Let $Z_{ij}(t)$ denote the number of class $i$ customers in service in pool $j$ at the beginning time period $t$. We also denote $U_{ij}(t)$ as the number of class $i$ customers assigned to pool $j$ for time period $t$. Note that $U_{ij}(t)$'s are determined by our scheduling policy.
Then the number of class $i$ departures from pool $j$ at the end of time period $t$, $R_{ij}(t)$, follows a Binomial distribution with parameter $Z_{ij}(t)+U_{ij}(t)$ and $\mu_{ij}$.
Let $X_i(t)$ denote the number of class $i$ customers waiting in queue at the beginning of period $t$. Then we have the following system dynamics,
\begin{equation}\label{queue_transition}
\begin{split} 
&X_i(t+1)=X_i(t)+A_i(t)-\sum_{j=1}^J U_{ij}(t), ~~ \forall i \in[I]\\ 
&Z_{ij}(t+1)=Z_{ij}(t)+U_{ij}(t)-R_{ij}(t), ~~ \forall i \in[I], j\in[J]. 
\end{split}
\end{equation}
%where $A_i(t) \sim \text{Poisson}(\theta_i)$, $R_{ij}(t) \big| \big(Z_{ij}(t),U_{ij}(t)\big)  \sim \text{Binomial}\big(Z_{ij}(t)+U_{ij}(t),\mu_{ij}\big)$.
The state of the system is $\bm{s}(t)=(X_i(t), Z_{ij}(t) : i \in[I], j\in[J])\in\mathbb{N}^{I\times(J+1)}$. The action is $\bm{a}(t)=(U_{ij}(t) : i \in[I], j\in[J])\in\mathbb{N}^{I\times J}$.

The routing policy needs to satisfy the following constraints
\begin{align} \label{queue_constraint1}
U_{ij}(t) \in \mathbb{N},\ \sum_{j=1}^J U_{ij}(t)\le X_i(t),\ \forall i \in[I], j\in[J],\ \forall t \ge 0.
\end{align}
i.e., we can not schedule more customers than there are waiting, and
\begin{align} \label{queue_constraint2}
\sum_{i=1}^I Z_{ij}(t)+U_{ij}(t)\le N_j,\ \forall  j\in[J],\ \forall t \ge 0,
\end{align}
i.e., the number of customers in service can not exceed the capacity. Note that  constraints \eqref{queue_constraint1}-\eqref{queue_constraint2}  are hard constraints, i.e, they need to be satisfied path-by-path. 

Each class $i$ customer waiting in queue incurs a holding cost of $h_i$ per period. There is also a one-shot routing cost of $r_{ij}$ for scheduling a class $i$ customer to a pool $j$ server. The overall cost for period $t$ is given by
\[c\big(\bm{s}(t),\bm{a}(t)\big)=\sum_{i=1}^I h_i X_i(t) +\sum_{i=1}^{I}\sum_{j=1}^{J} r_{ij}U_{ij}(t).\]
Our goal is to minimize the cumulative discounted costs: 
\[(1-\gamma)\cdot  \E^\pi \left[ \sum_{t=0}^\infty \gamma^t \cdot c(\bm{s}(t),\bm{a}(t))\right].\]

The problem we consider here is a weakly coupled MDP with $I$ sub-problems, where each sub-problem is an inverted-V model (i.e., a single customer class and multiple server pools).
In particular, for the $i$-th sub-problem, define state and action as $s^i(t)=(X_i(t), Z_{i1}(t),\ldots,Z_{iJ}(t))$ and $a^i(t)=(U_{i1}(t),\ldots,U_{iJ}(t))$. The transition dynamics of the $i$-th sub-system follows
\[
X_i(t+1)=X_i(t)+A_i(t)-\sum_{j=1}^J U_{ij}(t),  ~~ Z_{ij}(t+1)=Z_{ij}(t)+U_{ij}(t)-R_{ij}(t), ~~ \forall j\in[J].
\]
Given $s^i(t)$, the corresponding action space is defined as   
\[\cA^i(s^i(t))= \Big\{ \big\{ U_{ij}(t) \big \}_{j \in [J]} : U_{ij}(t) \in \mathbb{N},\ \sum_{j=1}^J U_{ij}(t) \le X_i(t),\  Z_{ij}(t)+U_{ij}(t) \le N_j,\ \forall j \in [J] \Big\}. \]
We also define the auxiliary cost function 
\[b^i(s^i(t),a^i(t))=\big(Z_{i1}(t)+U_{i1}(t),\ldots,Z_{iJ}(t)+U_{iJ}(t)\big)^\top  \in \mathbb{N}^J.\] 
Then the capacity constraints  \eqref{queue_constraint2} can be expressed as 
\begin{align*} 
\sum_{i=1}^I  b^i(s^i(t),a^i(t))  \le (N_1,\ldots,N_J)^\top,
\end{align*}
which takes the same form as the linking constraint in \eqref{joint_constraint}. 
%The dynamics of each subproblem is illustrated in Figure \ref{dynamics} (b).

%This problem is also referred to as ``skill-based routing problem''  in literature. It finds broad applications in service operations management. For example, \cite{Dai} use this model to study inpatients management. In hospital managements, inpatient beds are usually grouped into several wards (server pools), and each ward is assigned to serve patients (customers) from certain ``primary'' specialties. Without loss of generality, let  pool $i$ be the designated pool to inpatients of class $i$ with routing cost $r_{ii}=0, \forall i \in [I]$. In order to avoid the excessive admission delay, hospital managers may assign inpatient to an available non-primary bed if her designated  pool is full. However, such an ``overflow'' decision may be undesirable, since empirical evidence has shown that off-service placement can lead to adverse outcomes like a higher readmission rate \citep{Song}. This tradeoff is characterized by the nonzero routing costs $r_{ij}$ for $j\ne i$. \cite{Dai} apply approximate dynamic programming to solve the MDP directly. However, in contrast to our approach, they do not exploit the inherent coupling structure of this problem. 

There are three important features of the problem that we attempt to address in this section: 1) non-preemptive routing; 2) both class-and-pool dependent service rate; 3) routing cost (overflow cost). The first two features require us to keep track of a very high dimensional state space, i.e., $I(J+1)$. The third feature has not been extensively studied in the literature.

We next introduce two heuristic policies adapted from policies developed in the literature.
For PSS with multiple classes of customers and multiple pools of servers, a myopic policy called the $c\mu$-rule (or generalization of it), has been shown to be asymptotically optimal in some systems, where the goal is to minimize the holding cost \citep{mandelbaum2004scheduling}. The idea is to minimize the instantaneous cost-reduction rate at each decision epoch. 
Another policy is called the max-pressure policy, which is known to be throughput optimal and asymptotically cost optimal for some forms of convex holding cost \citep{stolyar2004maxweight,dai2008asymptotic}.
We next consider modified versions of the above routing policies, which take the routing costs into account \citep{queue_survey}. 
At each decision epoch $t$, we choose $U_{ij}(t)$'s that solve the following optimization problem:
\begin{align*}
\max_{U_{ij}(t)}& \quad   \sum_{i=1}^I\sum_{j=1}^J \omega_{ij}(t) U_{ij}(t) \\
\text{s.t.} &\quad \sum_{j=1}^J U_{ij}(t)\le X_i(t),\ \forall i \in[I], j\in[J],\ \forall t \ge 0, \\
& \quad   \sum_{i=1}^I Z_{ij}(t)+U_{ij}(t)\le N_j,\ \forall  j\in[J], \ \forall t \ge 0,\\
& \quad U_{ij}(t) \in \mathbb{N}\ \ \forall i \in[I], j\in[J],\ \forall t \ge 0,
\end{align*}
where $\omega_{ij}(t)$'s are some modified instantaneous costs we introduce next. We consider two different forms of $w_{ij}(t)$'s.
The first one sets $\omega_{ij}(t)=h_i-r_{ij}$, which is adapted from the $c\mu$-rule. We refer to this policy as the modified $c\mu$-rule. The second one sets $w_{ij}(t)=h_iX_i(t) - r_{ij}$, which is adapted from the max-pressure policy. We refer to this policy as the modified max-pressure policy.

%In $h\mu$-rule and max-pressure rule, the weights are defined as $\omega_{ij}=h_i\mu_{ij}$  and $\omega_{ij}=h_i\mu_{ij} X_i(t)$, respectively. Intuitively, the $h\mu$-rule gives high priority to customers with large holding cost, and aims to schedule them to pools with fast service rate.  The  max-pressure rule balances the tradeoff between holding cost and system congestion by incorporating the queue length into weights. Under zero routing cost $r_{ij}=0$, the $h\mu$-rule is proven to be optimal for single customer class problems. The max-pressure rule is asymptotic throughput optimal under heavy traffic regime \citep{mpp}. However, when there is nonzero routing cost, the extra cost structure complicates problems significantly. Few result about the optimality is known. Commonly-used heuristics modify the $h\mu$ and max-pressure rules by subtracting the routing costs in the corresponding weights.  Specifically, the modified $h\mu$-rule uses weights $\omega_{ij}=h_i\mu_{ij}-r_{ij}$ while the modified max-pressure rule uses weights $\omega_{ij}=h_i\mu_{ij}X_i(t)-r_{ij}$. \cite{queue_survey} compare the performance of those heuristics under different settings numerically. In later numerical studies, we mainly compare the performance of our policy with the modified $h\mu$ rule and the modified max-weight policy.  

\subsection{Solution method}
We consider the CMDP relaxation of the weakly coupled MDP:
\[\begin{split}
&\min_\pi\ (1-\gamma)\cdot  \E^\pi \Big[ \sum_{t=0}^\infty \gamma^t \cdot c(\bm{s}(t),\bm{a}(t))\Big]\\
\mbox{ s.t. } & (1-\gamma)\cdot \E\Big[ \sum_{i=1}^I\sum_{t=0}^{\infty}\gamma^t \cdot  d^i(s^i(t),a^i(t)) \Big] \le (N_1,\ldots,N_J)^\top,
\end{split}\] 
and apply the primal-dual algorithm to solve it. 
The decoupling allows us to translate the original problem to $I$ sub-problems. In particular, in each iteration, 
we use regularized policy iteration to update the scheduling policy for a single-class multi-pool system with modified instantaneous cost: 
\[c_{\lambda}^i\big(s^i(t),a^i(t)\big)=h_i X_i(t)+\sum_{j=1}^J r_{ij}U_{ij}(t)+\sum_{j=1}^J \lambda_j\big(Z_{ij}(t)+U_{ij}(t)\big)\]
for the $i$-th sub-problem.

Even with the decomposition, the state and policy spaces are still too large in this case. 
We next introduce some further approximations to reduce the dimension of the problem.
We shall omit the index $i$ in subsequent discussions as the development focuses on each sub-problem.

{\bf Policy space reduction:}
For each sub-problem, the policy space is still prohibitive. To see this, consider a system with $3$ pools and $30$ servers in each pool. When the queue length is $90$ and all pools are empty, there are roughly $30^3$ feasible actions. To overcome the challenge, we reduce the action space to only include priority rules. State-dependent extreme policies have been shown to be asymptotically optimal in the scheduling of PSS due to the linear system dynamics and linear holding costs \citep{harrison2004dynamic}. 
Denote $-1$ as the waiting option. The priority rule is denoted by a priority list that ends with $-1$. For example, priority $(1,2,-1)$ means pool 1 is preferred to pool 2, which is preferred to waiting. When following priority $(1,2,-1)$, we first assign as many customers to pool 1 as possible. If there are still customers waiting after pool 1 assignment, we start assigning them to pool 2. After that, if there are still customers waiting, we keep them in the queue. We denote this reduced policy space as $\tilde \cA$.

{\bf Value function approximation:} 
In our policy iteration step, given a policy $\pi$, we need to estimate the function $Q^{\pi,\lambda}(s,a)$ for all $s \in \cS$, $a\in \tilde \cA$, where the state $s=(x,z_1,\ldots,z_J)$.
Due to the large state space, we can not enumerate all the states to evaluate the value function. Instead, we use value function approximation with quadratic basis.
The idea is to find $\theta^{\pi,a} \in \mathbb{R}^{(J+1)^2+1}$ such that  
\[
Q^{\pi,\lambda}(s,a) \approx \langle \phi(s),  \theta^{\pi,a} \rangle.
\]
where $\phi(s)$ is the quadratic basis.
To obtain $\theta^{\pi,a}$ at each iteration, we first randomly sample $M$ states $\{s_i\}_{i \in [M]}$ and use Monte Carlo simulation to estimate $ Q^{\pi,\lambda}(s_i,a)$.  Then, set
\[\theta^{\pi,a}=\text{argmin}_\theta \Big\{ \frac{1}{M}\cdot \sum_{i=1}^M (Q^{\pi,\lambda}(s_i,a)-\langle \phi(s_i),\theta \rangle)^2\Big\}. 
\]

\subsection{Experiment Results}
For the numerical experiments, we consider a similar setting as that in \cite{Dai}, which is motivated by hospital inpatient-flow management.
In particular, we consider a network with 3 classes of customers and 3 pools of servers. Pool $i$ is considered the primary pool for class $i$ customers with $r_{ii}=0$, $\forall i \in [I]$. 
The major difference between our model and the model considered in \cite{Dai} is that we allow the service rates to vary for different server types, i.e, $\mu_{ij}$ depends on both $i$ and $j$. This captures the potential slowdown effect due to off-service placement \citep{Song}. %We also allow $r_{ij}$ to be different for different pool $j$. 
 
For the system parameters, we set the arrival rates $(\theta_1,\theta_2,\theta_3)=(12, 16, 20)$, the holding costs $(h_1,h_2,h_3)=(3,2,1)$,
the pool sizes $ (N_1,N_2,N_3)=(40,50,60) $, and the service rates 
\[
(\mu_{11},\mu_{12},\mu_{13})=(0.3,0.25,0.2), (\mu_{21},\mu_{22},\mu_{23})=(0.15,0.3,0.2), (\mu_{31},\mu_{32},\mu_{33})=(0.25,0.1,0.4).
\]
We run two sets of experiments, corresponding to large routing/overflow costs: 
\begin{equation}\label{eq:large_cost}
(r_{11},r_{12},r_{13})=(0,2,2), (r_{21},r_{22},r_{23})=(3,0,3), (r_{31},r_{32},r_{33})=(1,1,0),
\end{equation}
and small routing costs:
\begin{equation} \label{eq:small_cost}
(r_{11},r_{12},r_{13})=(0,0.2,0.2), (r_{21},r_{22},r_{23})=(0.3,0,0.3), (r_{31},r_{32},r_{33})=(0.1,0.1,0).
\end{equation}
Note that for class $i$ customers, the primary server pool $i$ has the largest service rate and zero routing cost. For customer class $i$, we define its nominal traffic intensity as $\rho_i=\theta_i/(N_i\mu_{ii})$.
Then the nominal traffic intensity of the three classes are
$\rho_1=1$, $ \rho_2=16/15$, and $\rho_3=5/6$. This indicates that the first two classes are unstable if we do not do any ``overflow". 
%We also set the routing costs comparable to the holding costs to make the tradeoff between waiting in queue and scheduling to non-primary pools more subtle. 

We initialize the system with $X_i(0)=50$ and $Z_{11}(0)=20$, $Z_{22}(0)=30$, $Z_{33}(0)=40$, and $Z_{ij}(0)=0$ for $i\neq j$, $i,j\in[3]$.
We compare the performance of our policy with the two benchmark policies for problems with different routing costs and discount rates. %, $\gamma=0.9,0.95,0.99$. 

When constructing the policy space for our primal-dual algorithms, %we restrict ourselves to priority policies. 
because each customer class has a primary server pool with the fastest service rate and zero routing cost, we always give the primary pool the highest priority. In particular, the action spaces for three classes are defined as 
%\[\mathcal{A}_1=\{ (1,-1),(1,2,-1),(1,3,-1),(1,2,3,-1),(1,3,2,-1)\}. \] 
%Recall that $(1,2,-1)$ means pool 1 is preferred to pool 2, which is preferred to waiting. 
%Similarly, the action spaces of classes 2 and 3 are defined as
\[\begin{split}
\mathcal{A}_1&=\{ (1,-1),(1,2,-1),(1,3,-1),(1,2,3,-1),(1,3,2,-1)\},\\
\mathcal{A}_2&=\{(2,-1),(2,1,-1),(2,3,-1),(2,1,3,-1),(2,3,1,-1) \}, \mbox{ and }\\
\mathcal{A}_3&=\{ (3,-1),(3,2,-1),(3,1,-1),(3,2,1,-1),(3,1,2,-1)\}, \mbox{ respectively.} 
\end{split}\] 

In our primal-dual update, we use the constant stepsize $0.1$.
When using simulation to estimate the value function, we truncate at $T=100, 150, 800$ for $\gamma=0.9,0.95,0.99 $ respectively. This ensures that $\gamma^T\approx 10^{-4}$, i.e., the truncation errors are almost negligible. 
When fitting the parameters for the quadratic value function approximation, we sample $1000$ states and use simulation to estimate the $Q$-function at these states.
For each value of $\gamma$, we start with the Lagrangian multipliers $\lambda_0=(10,10,10)$ and run the prima-dual algorithm for $30$ iterations, and take the policy obtained in the last iteration. Note that this policy may not be feasible to the original weakly coupled MDP. In order to obtain a feasible policy, we adopt the following modification. In each period, for each pool, when the number of scheduled customers exceeds the capacity, the primary customers are prioritized for admission. We then admit the ``overflowed'' customers uniformly at random until the capacity is reached. The customers who are not admitted to service will be sent back to their corresponding queues and wait for the next decision epoch. For example, suppose that there are $20$ servers available in pool 1 but the policy schedules $(15,5,5)$ customers from the three classes to this pool. The modified policy first admits the $15$ customers from class $1$ and then randomly picks $5$ among the $10$ customers of classes $2$ and $3$ to admit.  

Given a policy, to evaluate its performance, we estimate the cumulative discounted costs from $500$ independent replications of the system over $T$ periods of time. The results are summarized in Tables \ref{tab1} and \ref{tab2} . 

\begin{table}[!htbp]
\caption{Cumulative discounted costs under different policies with large routing costs \eqref{eq:large_cost} under different discount factors. (Numbers in bracket are the standard errors from simulation estimation.)}\label{tab1}
\centering
\begin{tabular}{c|c|c|c}
\hline
  & $\gamma=0.90$ & $\gamma=0.95$ & $\gamma=0.99$ \\
\hline
modified $c\mu$-rule &  $ \quad 270.49 \quad $    &  $\quad 286.11 \quad $ &  $ \quad 467.13  \quad $\\
& $ \quad (1.50) \quad $    &  $\quad (2.07) \quad $ &  $ \quad (4.26) \quad $\\
\hline
modified max-pressure rule & $271.62$  & $269.19$ & $278.31 $ \\
& $(1.25)$  & $(1.74)$ & $(2.31)$ \\
\hline
primal-dual algorithm & $\mathbf{252.77}$  & $\mathbf{243.87}$ & $\mathbf{218.95}$ \\
& $(1.47)$  & $(2.26)$ & $(2.11)$ \\
\hline
\end{tabular}
\end{table}

%\begin{table}[!htbp]
%\caption{Accumulated discounted costs under different policies with moderate routing cost under different discount factors. (Numbers in bracket are the standard errors from simulation estimation.)}\label{tab2}
%\centering
%\begin{tabular}{|c|c|c|c|}
%\hline
%  & $\gamma=0.90$ & $\gamma=0.95$ & $\gamma=0.99$ \\
%\hline
%modified $c\mu$-rule &  $ \quad 459.96 \quad $    &  $\quad 737.60 \quad $ &  $ \quad 3084.92 \quad $\\
%& $ \quad (0.88) \quad $    &  $\quad (1.20 ) \quad $ &  $ \quad (2.87 ) \quad $\\
%\hline
%modified max-pressure rule & $275.53$  & $283.34$ & $325.16$ \\
%& $(1.33)$  & $(1.62)$ & $(2.12)$ \\
%\hline
%primal-dual algorithm & $\mathbf{263.97}$  & $\mathbf{244.67}$ & $\mathbf{282.70}$ \\
%& $(1.48)$  & $(1.75)$ & $(3.38)$ \\
%\hline
%\end{tabular}
%\end{table}

\begin{table}[!htbp]
\caption{Cumulative discounted costs under different policies with small routing costs \eqref{eq:small_cost} under different discount factors. (Numbers in bracket are the standard errors from simulation estimation.)}\label{tab2}
\centering
\begin{tabular}{c|c|c|c}
\hline
  & $\gamma=0.90$ & $\gamma=0.95$ & $\gamma=0.99$ \\
\hline
modified $c\mu$-rule &  $ \quad 232.22 \quad $    &  $\quad 230.54 \quad $ &  $ \quad 266.81 \quad $\\
& $ \quad (1.20) \quad $    &  $\quad (1.71) \quad $ &  $ \quad (3.65) \quad $\\
\hline
modified max-pressure rule & $260.89$  & $266.83$ & $308.89$ \\
& $(1.21)$  & $(1.77)$ & $(2.06)$ \\
\hline
primal-dual algorithm & $\mathbf{253.53}$  & $\mathbf{251.20}$ & $\mathbf{210.34}$ \\
& $(1.50)$  & $(2.40)$ & $(2.24)$ \\
\hline
\end{tabular}
\end{table}

We observe that the policies obtained via the primal-dual algorithm performs  well. It outperforms the two benchmark policies in most cases.
When the routing cost is large (Table \ref{tab1}), the cost under the modified $c\mu$-rule increases substantially as the discount rate $\gamma$ increases. When taking a closer look at $w_{ij}(t)$'s, we note that in this case, $w_{21}(t)=-2.7$ and $w_{23}(t)=-2.6$. This implies that the modified $c\mu$-rule would never overflow class 2 customers. As a result, the system is unstable, i.e., the class $2$ queue blow up as $t$ increases. (The cumulative discounted cost is well-defined as the discount rate decays exponentially in $t$ while the queue length grows linearly in $t$.) 
The modified max-pressure is able to achieve reasonably good performance in this case. 
When $\gamma$ is small, our algorithm achieves comparable (slightly better) performance as the max-pressure policy. When $\gamma$ is large, i.e, $\gamma=0.99$,
our policy is able to achieve a substantially lower cost than the max-pressure policy, i.e., a 21\% cost reduction. This is because the max-pressure policy only starts overflowing when the queues are large enough. In this example where overflow is necessary to achieve system stability, we need more aggressive overflow. Our policy is able to ``learn" this through the primal-dual training.

When the overflow cost is small (Table (\ref{tab2}), the modified $c\mu$-rule is able to achieve better performance than the modified max-pressure policy. Note that in this case, all $w_{ij}(t)$'s are nonnegative for both the modified $c\mu$-rule and the modified max-pressure policy (when $ X_i(t)>0$). When $\gamma$ is small, our policy achieves comparable performance as the modified $c\mu$-rule, when $\gamma$ is large, i.e., $\gamma=0.99$, our policy can achieve a $21\%$ cost reduction over the modified $c\mu$-rule. This suggests that overflow needs to be exercised carefully.

%{\color{red} We observe that policies obtained via the primal-dual algorithm outperform the two benchmark policies. In addition, unlike the two benchmark policies, the performance of the policy obtained from the primal-dual algorithm improves as the discount factor decreases. Take the system with large routing costs for an example. 
%When taking a closer look at $w_{ij}(t)$, $j\neq i$, under the modified $c\mu$-rule and large routing cost, we note that they are all negative, e.g. $w_{12}(t)=3\times 0.25 -2=-1.75 $ and $w_{13}(t)=3\times 0.2 -2=-1.4$. This implies that the modified $c\mu$-rule would never overflow class 1 customers. 
%As a result, the system is not stable. Thus, the cumulative discounted cost increases as $\gamma$ increases. 
%This result also suggests that modified $c\mu$-rule may not be a good policy when having large routing cost. When  routing cost is small, the modified $c\mu$-rule gives positive indices for overflow decisions and makes system stable. However, our policies still improve system's performance when $\gamma$ is large. For the modified max-pressure rule, overflow only happens when the queue grows large enough. In contrast, the policies derived from our prima-dual algorithm do more aggressive overflow and hence, result in a lower cost in this case.

\begin{figure}[htbp]
\centering
\subfigure[Fix $Z_{12}=0, Z_{13}=0$, and vary $Z_{11}$]{
\begin{minipage}[t]{0.5\linewidth}
\centering
\includegraphics[width=3.3in]{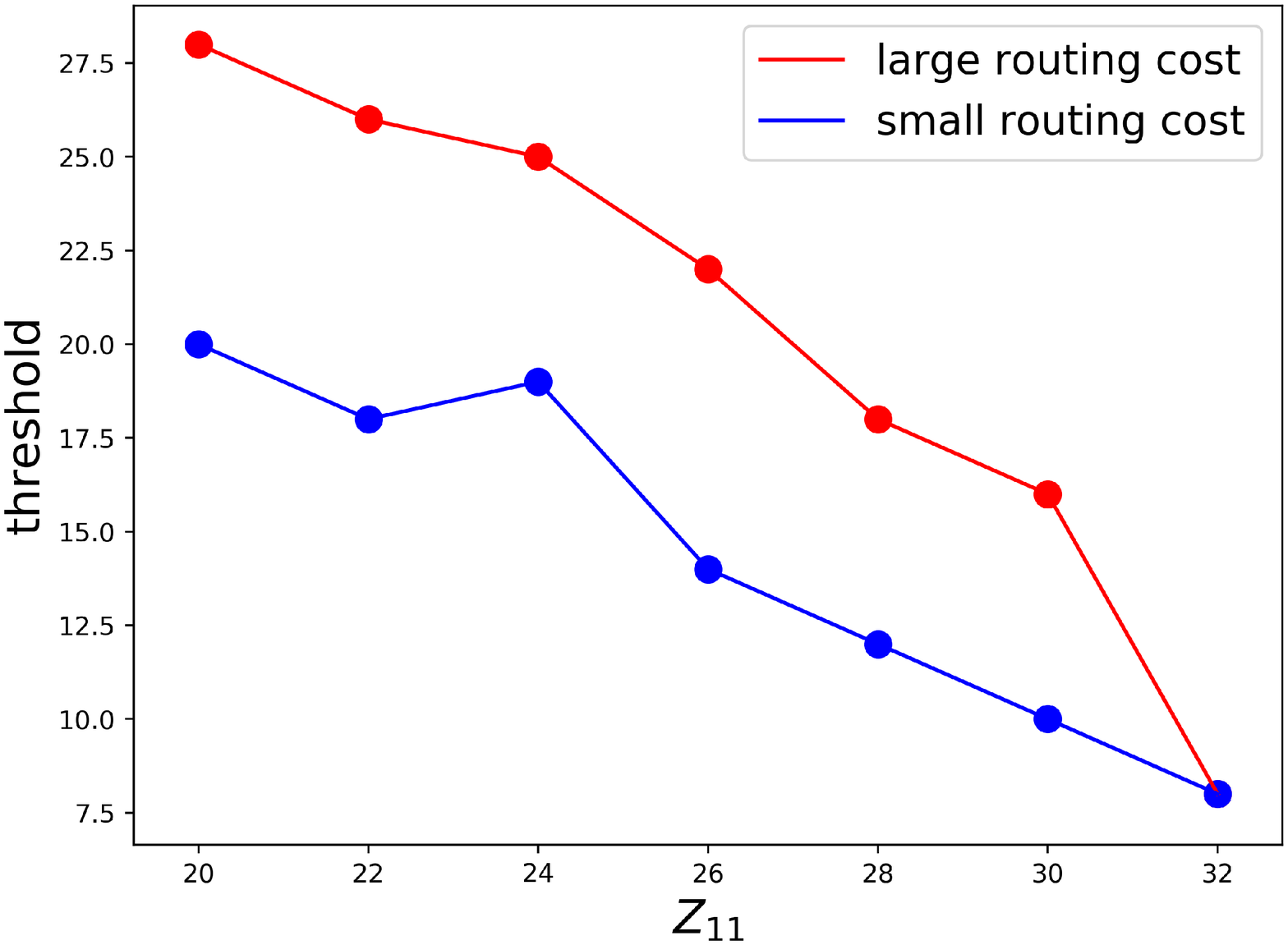}
%\caption{fig1}
\end{minipage}%
}%
\subfigure[Fix $Z_{21}=0, Z_{23}=0$, and vary $Z_{22}$]{
\begin{minipage}[t]{0.5\linewidth}
\centering
\includegraphics[width=3.3in]{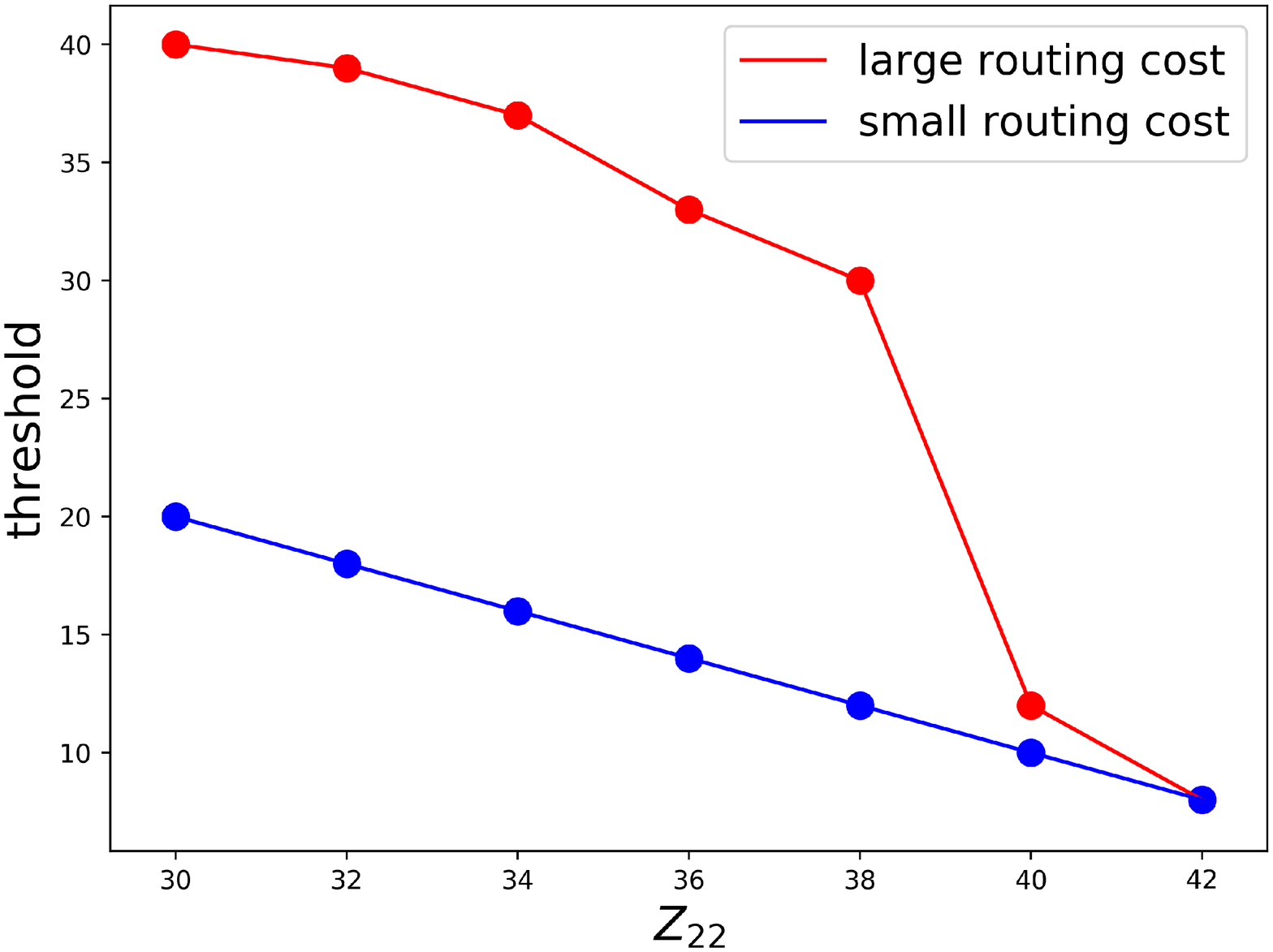}
%\caption{fig2}
\end{minipage}%
}%
\centering
\caption{The threshold for class 1 and 2 queues when starting overflowing.}  \label{thred}
\end{figure}

We next discuss the structure of the policies obtained via primal-dual algorithm. 
We observe that our policies in general follow a threshold structure: overflow customers only when the queue length exceeds some threshold. However, the thresholds are highly dependent on the states of the system. 
Take the scheduling policy for class $1$  and $2$ customers with discount rate $\gamma=0.9$ as an example. In Figure \ref{thred}, we plot the values of the threshold of starting overflowing for different values of $Z_{11}$'s and $Z_{22}$'s. We observe that holding $Z_{12} $ and $Z_{13}$ fixed, as $Z_{11}$ increases, the threshold for overflow decreases. Similarly, holding $Z_{21} $ and $Z_{23}$ fixed, as $Z_{22}$ increases, the threshold for overflow also decreases.

\section{Conclusion and Future Directions} \label{conclusion}
In this work, we propose a sampling-based primal-dual algorithm to solve CMDPs. Our approach alternatively applies regularized policy iteration to improve the policy and subgradient ascent to maintain the constraints. The algorithm achieves $O(\log(T)/\sqrt{T})$ convergence rate and only requires one policy update at each primal-dual iteration.
Our algorithm also enjoys the decomposability property for weakly coupled CMDPs.
We demonstrate the applications of our algorithm to solve two important operations management problems with weakly coupled structures: multi-product inventory management and multi-class queue scheduling.

In Section \ref{sec:queue}, we also show the good empirical performance of our algorithm to solve an important class of weakly coupled MDPs. 
This opens two directions for future research.
First, it is be important to quantify the optimality gap between the weakly coupled MDP and its CMDP relaxation theoretically. The gap can be large in some problems as demonstrate in \cite{Adelman}. It would be interesting to establish easy-to-verify conditions about when the gap is small. 
Second, the policy obtained via the Lagrangian relaxation may not satisfy the hard constraints in the original MDP. One approach to overcome the issue is to use more stringent thresholds when defining constraints in the CMDP relaxation \citep{Balseiro}. The other approach is to modify the CMDP based policies to construct good MDP policies. For example, \cite{brown2020index} study a dynamic assortment problem and propose an index heuristic from the relaxed problem, and show that the policy achieves asymptotic optimality. In Section \ref{sec:queue}, we apply a rather straightforward modification to the CMDP based policy in order to satisfy the hard constraints in the original MDP.
In general, how to ``translate" the policy derived based on the relaxed problem to the original MDP would be an interesting research direction. 

%Lastly, our algorithm can be viewed as an analog to the Lagrangian method in numerical optimization. It is well-known that the augmented Lagrangian method usually performs better (e.g. more stable and faster convergence rate). How to extend the idea of augmented Lagrangian to the CMDP setting is an interesting future research direction.

\bibliographystyle{informs2014}
\bibliography{reference}

\begin{thebibliography}{40}
\providecommand{\natexlab}[1]{#1}
\providecommand{\url}[1]{\texttt{#1}}
\providecommand{\urlprefix}{URL }

\bibitem[{Achiam et~al.(2017)Achiam, Held, Tamar, \protect\BIBand{}
  Abbeel}]{cpo}
Achiam J, Held D, Tamar A, Abbeel P (2017) Constrained policy optimization.
  \emph{arXiv preprint arXiv:1705.10528} .

\bibitem[{Adelman \protect\BIBand{} Mersereau(2008)}]{Adelman}
Adelman D, Mersereau AJ (2008) Relaxations of weakly coupled stochastic dynamic
  programs. \emph{Operations Research} 56(3):712--727.

\bibitem[{Altman(1999)}]{alt}
Altman E (1999) \emph{Constrained Markov decision processes}, volume~7 (CRC
  Press).

\bibitem[{Balseiro et~al.(2019)Balseiro, Brown, \protect\BIBand{}
  Chen}]{Balseiro}
Balseiro SR, Brown DB, Chen C (2019) Dynamic pricing of relocating resources in
  large networks. \emph{ACM SIGMETRICS Performance Evaluation Review}
  47(1):29--30.

\bibitem[{Bertsekas \protect\BIBand{} Scientific(2015)}]{bertsekas2015convex}
Bertsekas DP, Scientific A (2015) \emph{Convex optimization algorithms} (Athena
  Scientific Belmont).

\bibitem[{Bertsimas \protect\BIBand{} Orlin(1994)}]{BO}
Bertsimas D, Orlin JB (1994) A technique for speeding up the solution of the
  {L}agrangian dual. \emph{Mathematical Programming} 63(1-3):23--45.

\bibitem[{Bhatnagar \protect\BIBand{} Lakshmanan(2012)}]{Bhatnagar}
Bhatnagar S, Lakshmanan K (2012) An online actor--critic algorithm with
  function approximation for constrained Markov decision processes.
  \emph{Journal of Optimization Theory and Applications} 153(3):688--708.

\bibitem[{Borkar(2005)}]{Borkar}
Borkar VS (2005) An actor-critic algorithm for constrained Markov decision
  processes. \emph{Systems \& control letters} 54(3):207--213.

\bibitem[{Brown \protect\BIBand{} Smith(2020)}]{brown2020index}
Brown DB, Smith JE (2020) Index policies and performance bounds for dynamic
  selection problems. \emph{Management Science} .

\bibitem[{Bubeck(2014)}]{bubeck}
Bubeck S (2014) Convex optimization: Algorithms and complexity. \emph{arXiv
  preprint arXiv:1405.4980} .

\bibitem[{Caramanis et~al.(2014)Caramanis, Dimitrov, \protect\BIBand{}
  Morton}]{Caramanis}
Caramanis C, Dimitrov NB, Morton DP (2014) Efficient algorithms for
  budget-constrained Markov decision processes. \emph{IEEE Transactions on
  Automatic Control} 59(10):2813--2817.

\bibitem[{Chen et~al.(2020)Chen, Dong, \protect\BIBand{} Shi}]{queue_survey}
Chen J, Dong J, Shi P (2020) A survey on skill-based routing with applications
  to service operations management. \emph{Queueing Systems} 1--30.

\bibitem[{Chen \protect\BIBand{} Wang(2016)}]{SPD}
Chen Y, Wang M (2016) Stochastic primal-dual methods and sample complexity of
  reinforcement learning. \emph{arXiv preprint arXiv:1612.02516} .

\bibitem[{Chow et~al.(2018)Chow, Nachum, Duenez-Guzman, \protect\BIBand{}
  Ghavamzadeh}]{las}
Chow Y, Nachum O, Duenez-Guzman E, Ghavamzadeh M (2018) A {L}yapunov-based
  approach to safe reinforcement learning. \emph{Advances in neural information
  processing systems}, 8092--8101.

\bibitem[{Dai \protect\BIBand{} Shi(2019)}]{Dai}
Dai J, Shi P (2019) Inpatient overflow: An approximate dynamic programming
  approach. \emph{Manufacturing \& Service Operations Management}
  21(4):894--911.

\bibitem[{Dai et~al.(2008)Dai, Lin et~al.}]{dai2008asymptotic}
Dai JG, Lin W, et~al. (2008) Asymptotic optimality of maximum pressure policies
  in stochastic processing networks. \emph{The Annals of Applied Probability}
  18(6):2239--2299.

\bibitem[{Gattami(2019)}]{Gattami}
Gattami A (2019) Reinforcement learning for multi-objective and constrained
  Markov decision processes. \emph{arXiv preprint arXiv:1901.08978} .

\bibitem[{Geist et~al.(2019)Geist, Scherrer, \protect\BIBand{} Pietquin}]{rmdp}
Geist M, Scherrer B, Pietquin O (2019) A theory of regularized Markov decision
  processes. \emph{arXiv preprint arXiv:1901.11275} .

\bibitem[{Haarnoja et~al.(2017)Haarnoja, Tang, Abbeel, \protect\BIBand{}
  Levine}]{softq}
Haarnoja T, Tang H, Abbeel P, Levine S (2017) Reinforcement learning with deep
  energy-based policies. \emph{arXiv preprint arXiv:1702.08165} .

\bibitem[{Harrison \protect\BIBand{} Zeevi(2004)}]{harrison2004dynamic}
Harrison JM, Zeevi A (2004) Dynamic scheduling of a multiclass queue in the
  halfin-whitt heavy traffic regime. \emph{Operations Research} 52(2):243--257.

\bibitem[{Kakade \protect\BIBand{} Langford(2002)}]{kak}
Kakade S, Langford J (2002) Approximately optimal approximate reinforcement
  learning. \emph{ICML}, volume~2, 267--274.

\bibitem[{Le et~al.(2019)Le, Voloshin, \protect\BIBand{} Yue}]{bpl}
Le HM, Voloshin C, Yue Y (2019) Batch policy learning under constraints.
  \emph{arXiv preprint arXiv:1903.08738} .

\bibitem[{Liu et~al.(2019{\natexlab{a}})Liu, Cai, Yang, \protect\BIBand{}
  Wang}]{Liu}
Liu B, Cai Q, Yang Z, Wang Z (2019{\natexlab{a}}) Neural proximal/trust region
  policy optimization attains globally optimal policy. \emph{arXiv preprint
  arXiv:1906.10306} .

\bibitem[{Liu et~al.(2019{\natexlab{b}})Liu, Ding, \protect\BIBand{} Liu}]{IPO}
Liu Y, Ding J, Liu X (2019{\natexlab{b}}) Ipo: Interior-point policy
  optimization under constraints. \emph{arXiv preprint arXiv:1910.09615} .

\bibitem[{Mandelbaum \protect\BIBand{}
  Stolyar(2004)}]{mandelbaum2004scheduling}
Mandelbaum A, Stolyar AL (2004) Scheduling flexible servers with convex delay
  costs: Heavy-traffic optimality of the generalized c$\mu$-rule.
  \emph{Operations Research} 52(6):836--855.

\bibitem[{Miryoosefi et~al.(2019)Miryoosefi, Brantley, Daume~III, Dudik,
  \protect\BIBand{} Schapire}]{rlcc}
Miryoosefi S, Brantley K, Daume~III H, Dudik M, Schapire RE (2019)
  Reinforcement learning with convex constraints. \emph{Advances in Neural
  Information Processing Systems}, 14093--14102.

\bibitem[{Nedi{\'c} \protect\BIBand{} Ozdaglar(2009)}]{ssp}
Nedi{\'c} A, Ozdaglar A (2009) Subgradient methods for saddle-point problems.
  \emph{Journal of optimization theory and applications} 142(1):205--228.

\bibitem[{Neely(2011)}]{OFP}
Neely MJ (2011) Online fractional programming for {M}arkov decision systems.
  \emph{2011 49th Annual Allerton Conference on Communication, Control, and
  Computing (Allerton)}, 353--360 (IEEE).

\bibitem[{Nemirovski(2012)}]{md}
Nemirovski A (2012) Tutorial: Mirror descent algorithms for large-scale
  deterministic and stochastic convex optimization. \emph{Conference on
  Learning Theory (COLT)}.

\bibitem[{Schulman et~al.(2015)Schulman, Levine, Abbeel, Jordan,
  \protect\BIBand{} Moritz}]{trpo}
Schulman J, Levine S, Abbeel P, Jordan M, Moritz P (2015) Trust region policy
  optimization. \emph{International conference on machine learning},
  1889--1897.

\bibitem[{Schulman et~al.(2017)Schulman, Wolski, Dhariwal, Radford,
  \protect\BIBand{} Klimov}]{ppo}
Schulman J, Wolski F, Dhariwal P, Radford A, Klimov O (2017) Proximal policy
  optimization algorithms. \emph{arXiv preprint arXiv:1707.06347} .

\bibitem[{Shani et~al.(2019)Shani, Efroni, \protect\BIBand{} Mannor}]{Shani}
Shani L, Efroni Y, Mannor S (2019) Adaptive trust region policy optimization:
  Global convergence and faster rates for regularized mdps. \emph{arXiv
  preprint arXiv:1909.02769} .

\bibitem[{Singh \protect\BIBand{} Cohn(1998)}]{Singh}
Singh SP, Cohn D (1998) How to dynamically merge Markov decision processes.
  \emph{Advances in neural information processing systems}, 1057--1063.

\bibitem[{Sion et~al.(1958)}]{Sion}
Sion M, et~al. (1958) On general minimax theorems. \emph{Pacific Journal of
  mathematics} 8(1):171--176.

\bibitem[{Song et~al.(2020)Song, Tucker, Graue, Moravick, \protect\BIBand{}
  Yang}]{Song}
Song H, Tucker AL, Graue R, Moravick S, Yang JJ (2020) Capacity pooling in
  hospitals: The hidden consequences of off-service placement. \emph{Management
  Science} 66(9):3825--3842.

\bibitem[{Stolyar et~al.(2004)}]{stolyar2004maxweight}
Stolyar AL, et~al. (2004) Maxweight scheduling in a generalized switch: State
  space collapse and workload minimization in heavy traffic. \emph{The Annals
  of Applied Probability} 14(1):1--53.

\bibitem[{Sutton \protect\BIBand{} Barto(2018)}]{Sutton}
Sutton RS, Barto AG (2018) \emph{Reinforcement learning: An introduction} (MIT
  press).

\bibitem[{Tessler et~al.(2018)Tessler, Mankowitz, \protect\BIBand{}
  Mannor}]{rcpo}
Tessler C, Mankowitz DJ, Mannor S (2018) Reward constrained policy
  optimization. \emph{arXiv preprint arXiv:1805.11074} .

\bibitem[{Turken et~al.(2012)Turken, Tan, Vakharia, Wang, Wang,
  \protect\BIBand{} Yenipazarli}]{MPNP}
Turken N, Tan Y, Vakharia AJ, Wang L, Wang R, Yenipazarli A (2012) The
  multi-product newsvendor problem: Review, extensions, and directions for
  future research. \emph{Handbook of newsvendor problems}, 3--39 (Springer).

\bibitem[{Wang et~al.(2019)Wang, Cai, Yang, \protect\BIBand{}
  Wang}]{wang2019neural}
Wang L, Cai Q, Yang Z, Wang Z (2019) Neural policy gradient methods: Global
  optimality and rates of convergence. \emph{arXiv preprint arXiv:1909.01150} .

\end{thebibliography}

\newpage
\begin{APPENDIX}{Proof of Main Results}   \label{proof} 

%\subsection{Proof of Main Result}

The proof of Theorem \ref{thm1} relies the following lemma, which upper and lower bounds the movement of the Lagrangian after a single iteration/update of the policy and the Lagrangian multipliers. 
\begin{lemma} \label{descent_lemma}
  Let $\{(\pi_m,\lambda_m)\}_{m\ge 0}$ be the sequences of stationary policies and Lagrangian multipliers generated by Algorithm \ref{main_alg}. Then for arbitrary $\lambda \in \R^K_+$ and $\pi \in \Pi_S$, we have the upper bound
  \begin{align*}
   L(\pi_m,\lambda)-L(\pi_m,\lambda_m) \le (2\eta_{m})^{-1}\cdot \big(\|\lambda-\lambda_m \|^2-\|\lambda-\lambda_{m+1} \|^2\big) +  \eta_{m}/2 \cdot \big \|  \partial_\lambda L(\pi_m,\lambda_{m}) \big \|^2, 
  \end{align*}
  and the lower bound 
  \begin{align*}
  & L(\pi,\lambda_{m})-L(\pi_m,\lambda_m)\\
  & \qquad  \ge \big((1-\gamma)\eta_{m}\big)^{-1} \Big(\Phi^{\pi}(\pi\|\pi_{m+1})-\Phi^{\pi}(\pi\|\pi_m)  \Big)-\frac{\eta_{m}}{8(1-\gamma)}\cdot \big(\sup_{s\in \cS} \sup_{a\in \cA} |Q^{\lambda_m,\pi_m}(s,a) |\big)^2.
  \end{align*}  
\end{lemma}

Before we prove Lemma \ref{descent_lemma}, we first present two auxiliary lemmas.
The first lemma (Lemma \ref{general_MD}) is rather standard. A similar version of the result can be found in Proposition 3.2.2 in \cite{bertsekas2015convex}.  For self-completeness, we still provide the proof here.  

\begin{lemma} \label{general_MD}
  Let $f$ be a proper convex function on a space $\Omega$ (not necessary a Euclidean space). Let $\cC $ be an open set in $\Omega$, and $\Psi_{\xi}(\cdot\|\cdot)$ be the Bregman divergence induced by a strictly convex function $\xi$ on $\Omega$. For an arbitrary constant $\eta>0$ and a point $x_0 \in \Omega$, define 
  $$ x^*=\argmin_{x\in \cC} \Big \{f(x)+\frac{1}{\eta} \Psi_\xi\big(x\|x_0\big) \Big\}.$$ 
  Then we have 
  $$
  f(x)-f(x^*) \ge \frac{1}{\eta}  \Big(\Psi_\xi\big(x^*\|x_0\big)+\Psi_\xi\big(x\|x^*\big)-\Psi_\xi\big(x\|x_0\big) \Big),\ \forall\ x\in \Omega.
  $$
  By symmetry, for a concave function $g$ on $\Omega$ and 
  $$\hat x^*=\argmax_{x\in \cC} \Big \{g(x)-\frac{1}{\eta} \Psi_\xi\big(x\|x_0\big) \Big\}.$$ 
  Then 
  $$
 g(x)-g(\hat x^*) \le -\frac{1}{\eta}  \Big(\Psi_\xi\big(\hat x^*\|x_0\big)+\Psi_\xi\big(x\| \hat x^*\big)-\Psi_\xi\big(x\|x_0\big) \Big),\ \forall\ x\in \Omega.
  $$
\end{lemma}

\proof{Proof of Lemma  \ref{general_MD}}
We first consider the minimization problem. Since $x^*$ minimizes the objective $f(x)+ \eta^{-1}\cdot  \Psi_{\xi}(x\|x_0)$ on set $\mathcal{C} $, there exists a subgradient of the form
$$
p^*=q^*+ \eta^{-1} \cdot \partial_x \Psi_{\xi}(x^*\|x_0)=q^*+\eta^{-1}\cdot \big( \nabla \xi(x^*)-\nabla \xi(x_0) \big)
$$  
such that 
$$
\langle p^*, x-x^* \rangle \ge 0, \ \forall x\in \mathcal{C}.
$$
Here $q^* \in \partial_x f(x^*) $ is some subgradient of $f(x)$ at $x^*$. 
As a result, by the property of subgradient, for all $x \in \mathcal{C}$, we have 
\begin{align*}
f(x) &  \ge f(x^*) +\langle q^*, x-x^*  \rangle \\
&\ge f(x^*)+ \eta^{-1}\cdot  \langle  \nabla \xi(x_0)-\nabla \xi(x^*), x-x^* \rangle  \\
&= f(x^*)+ \eta^{-1}\cdot   \Big(\Psi_\xi\big(x^*\|x_0\big)+\Psi_\xi\big(x\|x^*\big)-\Psi_\xi\big(x\|x_0\big) \Big),
\end{align*}
where the last equality follows from the definition of Bregman divergence, i.e., 
$$
\Psi_{\xi}(x\|y)=\xi(x)-\xi(y)- \big \langle \nabla \xi(y),  x-y \big \rangle. 
$$
For the maximization problem, we only need to consider $-g$ and apply above result. 
\Halmos
\endproof

 The next lemma is Lemma 6.1 in \citep{kak}. Given two policies, it characterizes the difference of expected accumulated costs as the inner product of the advantage function of one policy and the occupation measure of another policy. Note that the value function $V^{\pi}$ and the action-value function $Q^\pi$  of an MDP under policy $\pi$ are defined in \eqref{def:value_func} and \eqref{Q_function}.
\begin{lemma} \label{MDP_difference}
%  Consider an MDP with initial distribution $\mu_0$. Let $V^{\pi}(s)$, $Q^\pi(s,a)$, and $\nu^\pi$ be the value function, action-value function, and occupation measure associated with a stationary policy $\pi \in \Pi_S$, respectively. Then for 
For arbitrary policies $\pi,\pi' \in \Pi_S$, 
  $$
  \E_{s\sim \mu_0}\big[V^\pi(s)\big]- \E_{s\sim \mu_0}\big[V^{\pi'}(s)\big]=\frac{1}{1-\gamma} \E_{(s,a)\sim \nu^{\pi'}}\big[ Q^\pi(s,a)-V^\pi(s)\big].
  $$
where $\nu^{\pi'}(\cdot,\cdot) $ is the  occupation measure associated with $\pi'$.
\end{lemma}
%Lemma \ref{general_MD} is a standard result in the analysis of mirror descent algorithm. See for example, \cite{cvxo}, for more details. To prove Lemma \ref{descent_lemma}, we first rewrite the iteration in \eqref{iteration} as the argmax and argmin form in Lemma \ref{general_MD}. Then the results hold naturally with a tailored choice of the Bregman divergence $\Psi_{\xi}$. Lemma \ref{MDP_difference} is from \citep{kak}. It expresses the difference of accumulated costs as the inner product of the advantage function of one policy and the occupation measure of another policy, which enables us to use Cauchy-Schwartz inequality to separate those two parts. 

\proof{Proof of Lemma  \ref{descent_lemma}}
For the upper bound, note that because $L(\pi_m,\lambda) $ is linear in $ \lambda$,
$$
\lambda_{m+1}=\text{Proj}_{\Lambda_M}\big\{ \lambda_m+\eta_{m} \cdot\partial_\lambda L(\pi_m,\lambda_m)\big\}
$$
is equivalent to 
\begin{align*}
\lambda_{m+1}=\argmax_{\lambda \in {\Lambda_M}} \Big\{  L(\pi_m,\lambda)-\frac{1}{2\eta_{m}} \|\lambda-\lambda_m\|^2  \Big\}.
\end{align*}
Then, by Lemma \ref{general_MD}, we have 
\begin{align*}
  L(\pi_m,\lambda)-L(\pi_m,\lambda_{m+1}) & \le (2\eta_{m})^{-1}\big(\|\lambda-\lambda_m \|^2-\|\lambda-\lambda_{m+1} \|^2-\|\lambda_{m+1}-\lambda_m \|^2\big) \notag\\
  & \le  (2\eta_{m})^{-1} \big(\|\lambda-\lambda_m \|^2-\|\lambda-\lambda_{m+1} \|^2\big).
\end{align*} 
Next,
\begin{align*}
 L(\pi_m,\lambda)-L(\pi_m,\lambda_{m})
 \le& (2\eta_{m})^{-1} \big(\|\lambda-\lambda_m \|^2-\|\lambda-\lambda_{m+1} \|^2\big) +  L(\pi_m,\lambda_{m+1})-L(\pi_m,\lambda_{m})    \notag\\
  =&  (2\eta_{m})^{-1} \big(\|\lambda-\lambda_m \|^2-\|\lambda-\lambda_{m+1} \|^2\big) +  \big \langle  \partial_\lambda L(\pi_m,\lambda_{m})  ,\lambda_{m+1}-\lambda_m \big \rangle  \notag \\
  \le& (2\eta_{m})^{-1} \big(\|\lambda-\lambda_m \|^2-\|\lambda-\lambda_{m+1} \|^2\big) +  \eta_{m}/2 \cdot \big \|  \partial_\lambda L(\pi_m,\lambda_{m}) \big \|^2,
\end{align*} 
where the last inequality follows from the definition of $\lambda_{m+1}$ and the non-expansive property of the projection. Then we obtain the upper bound.

For the lower bound, recall that we update $\pi_m$ via  
\begin{align*}
  \pi_{m+1}(\cdot|s)=\argmin_{\pi(\cdot|s) \in \Delta_{\cA} } \Big \{ \big \langle Q^{\pi_m,\lambda_m}(s,\cdot), \pi(\cdot|s) \big\rangle+\frac{1}{\eta_{m}} \KL\big( \pi(\cdot|s) \| \pi_m(\cdot|s)\big) \Big \},
\end{align*}
for each state $s\in \cS$.
Then, for an arbitrary stationary policy $\pi' \in \Pi_S$,  we have 
\begin{align*}
  \pi_{m+1}=\argmin_{\pi \in \Pi_S } \Big \{    \E_{s\sim \nu_s^{\pi'}}\Big[  \big \langle Q^{\pi_m,\lambda_m}(s,\cdot), \pi(\cdot|s) \big\rangle+\frac{1}{\eta_{m}} \KL\big( \pi(\cdot|s) \| \pi_m(\cdot|s)\big) \Big] \Big \}
\end{align*}
where $\nu_s^{\pi'} $ is the state occupation measure associated with $\pi'$. %In Lemma \ref{general_MD}, 

Note that the space of the stationary policy, $\Pi_S $, can be represented as the product space of simplex $\Delta_{\cA}$.
Consider $\Omega:=\Pi_S= \big(\Delta_{\cA}\big)^{\otimes|\cS|}$ and let
\begin{align*}
  g(\pi)&:=  \E_{s\sim \nu_s^{\pi'}}\Big[  \big \langle Q^{\pi_m,\lambda_m}(s,\cdot), \pi(\cdot|s) \big\rangle \Big], \\
  \Psi_\xi(\pi)&:=\E_{s\sim \nu_s^{\pi'}}\Big[\KL\big(\pi(\cdot|s)\|\pi_m(\cdot|s)\big) \Big]=\Phi^{\pi'}(\pi\|\pi_m).
\end{align*}  
where $\Phi^{\pi'}$ is defined in \eqref{w-KL}.  Since $g(\pi)$ is linear in $\pi$, setting $\pi=\pi'$, by Lemma \ref{general_MD}, we obtain
\begin{align*}
  \E_{s\sim \nu_s^{\pi'}}\Big[  \big \langle Q^{\pi_m,\lambda_m}(s,\cdot), \pi'(\cdot|s) - \pi_{m+1}(\cdot|s) \big\rangle\Big] \ge \eta_{m}^{-1}  \Big (\Phi^{\pi'}(\pi_{m+1}\|\pi_m)+\Phi^{\pi'}(\pi'\|\pi_{m+1})-\Phi^{\pi'}(\pi'\|\pi_m) \Big),
\end{align*}
which can be equivalently written as 
\begin{align} \label{ineq01}
 &  \eta_{m}^{-1} \cdot \Big( \Phi^{\pi'}(\pi'\|\pi_{m+1}) - \Phi^{\pi'}(\pi'\|\pi_m) + \Phi^{\pi'}(\pi_{m+1}\|\pi_m)  \Big) \notag \\
 & \le \E_{s\sim \nu_s^{\pi'}}\Big[  \big \langle Q^{\pi_m,\lambda_m}(s,\cdot), \pi'(\cdot|s) - \pi_{m}(\cdot|s) \big\rangle\Big]+  \E_{s\sim \nu_s^{\pi'}}\Big[  \big \langle Q^{\pi_m,\lambda_m}(s,\cdot), \pi_m(\cdot|s) - \pi_{m+1}(\cdot|s) \big\rangle\Big].
\end{align} 

We next derive an upper bound for the right-hand side of inequalities \eqref{ineq01}. Let $\|\cdot \|_{\text{TV}}$ denotes the total variation norm of probability distributions.
First, for each state $s \in \cS$,  
\begin{align*} 
& \eta_{m} \cdot  \big \langle Q^{\pi_m,\lambda_m}(s,\cdot), \pi_m(\cdot|s) - \pi_{m+1}(\cdot|s) \big\rangle \\
\le&  \eta_{m}  \cdot \sup_{a\in \cA} \big|Q^{\pi_m,\lambda_m}(s,a) \big|\cdot \big\| \pi_m(\cdot|s) - \pi_{m+1}(\cdot|s) \big\|_{\text{TV}} \\ %\mbox{\ (by Cauchy-Schwartz inequality)}\\
\le&  \frac{\eta_{m}^2}{8}  \cdot \Big(\sup_{s\in \cS} \sup_{a\in \cA} \big|Q^{\pi_m,\lambda_m}(s,a) \big|\Big)^2+ 2\cdot  \big\| \pi_{m+1}(\cdot|s) - \pi_{m}(\cdot|s) \big\|^2_{\text{TV}} \\
\le&  \frac{\eta_{m}^2}{8}  \cdot \Big(\sup_{s\in \cS} \sup_{a\in \cA} \big|Q^{\pi_m,\lambda_m}(s,a) \big|\Big)^2+\KL\big( \pi_{m+1}(\cdot|s) \| \pi_{m}(\cdot|s) \big) \mbox{\ (by Pinsker's inequality)}.
\end{align*}
Hence, by taking the average, we obtain
\begin{equation} \label{ineq02}
\begin{split}
&\E_{s\sim \nu_s^{\pi'}}\Big[  \big \langle Q^{\pi_m,\lambda_m}(s,\cdot), \pi_m(\cdot|s) - \pi_{m+1}(\cdot|s) \big\rangle\Big]\\
\le& \frac{\eta_{m}}{8} \cdot \Big(\sup_{s\in \cS} \sup_{a\in \cA} \big|Q^{\pi_m,\lambda_m}(s,a) \big|\Big)^2+ \eta^{-1}_{m}\cdot  \Phi^{\pi'}(\pi_{m+1}\|\pi_m).
\end{split}
\end{equation}
Second, recall that $\nu^\pi(s,a)=\nu_s^\pi(s)\cdot \pi(a|s)$ and $V^\pi(s)=\langle Q^\pi(s,\cdot),\pi(\cdot|s) \rangle $. Then, by Lemma \ref{MDP_difference}, for the modified unconstrained MDP, we have 
\begin{align}  \label{ineq03}
 \E_{s\sim \nu_s^{\pi'}}\Big[  \big \langle Q^{\pi_m,\lambda_m}(s,\cdot), \pi'(\cdot|s) - \pi_{m}(\cdot|s) \big\rangle\Big] & =\E_{(s,a)\sim \nu^{\pi'}}\big[  Q^{\pi_m,\lambda_m}(s,a) - V^{\pi_m,\lambda_m}(s) \big]  \notag \\
 &=(1-\gamma) \cdot \big(L(\pi',\lambda_m)-L(\pi_m,\lambda_m)\big).
\end{align}
Finally, combining \eqref{ineq01}-\eqref{ineq03}, we obtain
\begin{align*}
 & L(\pi',\lambda_m)-L(\pi_m,\lambda_m) \\
\ge& \big((1-\gamma)\eta_{m}\big)^{-1} \Big(\Phi^{\pi'}(\pi'\|\pi_{m+1})-\Phi^{\pi'}(\pi'\|\pi_m) \Big) -\frac{\eta_{m}}{8(1-\gamma)} \big(\sup_{s\in \cS} \sup_{a\in \cA} |Q^{\pi_m,\lambda_m}(s,a) |\big)^2 .
\end{align*}
\Halmos 
\endproof

We are now ready to prove Theorem \ref{thm1}.
\proof{Proof of Theorem \ref{thm1}}
{\bf We prove the bound for $D(\bar \pi_T)-q$ first.} For this, we only need to establish an upper bound for
$\big \| [\partial_{\lambda}L(\bar{\pi}_T,\lambda) ]^+ \big\|.$ 

%Recall that we update the Lagrangian multipliers via 
%$$
%\lambda_{m+1}=\text{Proj}_{\Lambda_M}\Big\{ \lambda_m+\eta_{m} \cdot \partial_{\lambda}L(\pi_m,\lambda_m)\Big\}.
%$$
Since $L(\pi_m,\lambda)$ is linear in $\lambda$, we have 
\begin{align}  \label{001}
  L(\pi_m,\lambda_m)-L(\pi_m,\lambda^*)=(\lambda_m-\lambda^*)^\top \partial_{\lambda} L(\pi_m,\lambda_m).
\end{align} 
By the first part of Lemma \ref{descent_lemma}, for any $\lambda$, we have
%By the non-expansive property of the projection to a closed convex set and Assumption \ref{ass2}, for any $\lambda \in \Lambda_M$, we have 
%\begin{align*}
%  \|\lambda_{m+1}-\lambda\|^2 & \le \big\| \lambda_m+\eta_{m} \cdot \partial_{\lambda}L(\pi_m,\lambda_m)-\lambda \big\|^2 \\
%  & \le \|\lambda_{m}-\lambda\|^2+2\eta_{m+1}\cdot (\lambda_{m}-\lambda)^\top\partial_{\lambda}L(\pi_m,\lambda_m)+\eta_{m}^2 G^2,
%\end{align*}
%which can be equivalently expressed as
\begin{align}  \label{002}
\eta_{m}\cdot (\lambda-\lambda_{m})^\top\partial_{\lambda}L(\pi_m,\lambda_m) &= \eta_{m}\cdot ( L(\pi_m,\lambda)-L(\pi_m,\lambda_m))  \nonumber \\
 & \le \big( \|\lambda_{m}-\lambda\|^2-\|\lambda_{m+1}-\lambda\|^2  \big)/2+\eta_{m}^2 G^2/2.
\end{align}
On the other hand, by the saddle point property of $(\pi^*,\lambda^*)$, we also have  
\begin{align}  \label{003}
  L(\pi_m,\lambda^*) \ge L(\pi^*,\lambda^*).%=L^*_{\mu_0}.
\end{align} 
In the following, we denote by $L^*:=L(\pi^*,\lambda^*)$. By combining inequalities \eqref{001}-\eqref{003}, we obtain 
\begin{align*}
& \eta_{m} \cdot  (\lambda-\lambda^*)^\top \partial_{\lambda} L(\pi_m,\lambda_m)\\
=&\eta_{m}\cdot(\lambda-\lambda_m)^\top \partial_{\lambda} L(\pi_m,\lambda_m)+\eta_{m}\cdot (\lambda_m-\lambda^*)^\top \partial_{\lambda} L(\pi_m,\lambda_m) \\
 \le& \big( \|\lambda_{m}-\lambda\|^2-\|\lambda_{m+1}-\lambda\|^2  \big)/2+\eta_{m}^2 G^2/2 + \eta_{m} \cdot \big (L(\pi_m,\lambda_m)-L^*\big).
\end{align*}
By taking the telescope sum of above inequality, for any $\lambda \in \Lambda_M$, we have
\begin{equation} \label{004}
\begin{split}
  &\sum_{m=0}^{T-1} \eta_m \cdot  (\lambda-\lambda^*)^\top \partial_{\lambda} L(\pi_m,\lambda_m)\\
  \le& \big( \|\lambda_{0}-\lambda\|^2-\|\lambda_{T}-\lambda\|^2  \big)/2+ \Big(\sum_{m=0}^{T-1}\eta_m^2 /2\Big)\cdot G^2\\
  &+\sum_{m=0}^{T-1} \eta_m \cdot \big(L(\pi_m,\lambda_m)-L^*\big).
\end{split}
\end{equation}

For the left hand side of \eqref{004}, let
$$
\zeta_T:=\sum_{m=0}^{T-1} \eta_m \cdot \partial_{\lambda} L(\pi_m,\lambda_m)=\Big(\sum_{m=0}^{T-1} \eta_m\Big) \cdot \partial_\lambda L(\bar{\pi}_T,\lambda),
$$
where the last equality follows from the definition of $\bar{\pi}_T$ and the linearity of value function under the mixing operation. If $[\zeta_T]^+=0$, then the upper bound holds trivially. Otherwise, let 
$$
\tilde{\lambda}=\lambda^*+r\cdot \frac{[\zeta_T]^+}{\big\| [\zeta_T]^+ \big\|},
$$
where $r$ is the slackness constant in the definition of $\Lambda_M$ in \eqref{Lambda_area}. Then it is easy to see that $\tilde{\lambda} \in \Lambda_M$. By \eqref{004}, we have
\begin{align*}
(\tilde{\lambda}-\lambda^*)^\top \zeta_T \le \max_{\lambda \in \Lambda_M} \|\lambda-\lambda_0\|^2 /2 +\Big(\sum_{m=0}^{T-1} \eta_m^2 /2\Big)\cdot G^2 + \sum_{m=0}^{T-1} \eta_m \cdot \big(L(\pi_m,\lambda_m)-L^*\big).
\end{align*}
By the definition of $\tilde{\lambda}$, we also have 
$$
(\tilde{\lambda}-\lambda^*)^\top \zeta_T= r\cdot  \frac{([\zeta_T]^+)^\top\zeta_T}{\big\|[\zeta_T]^+ \big\|}=r\cdot \big\|[\zeta_T]^+ \big\| = r\cdot \Big(\sum_{m=0}^{T-1} \eta_m\Big) \cdot \big\|[ \partial_\lambda L(\bar{\pi}_T,\lambda) ]^+ \big\|.
$$
Hence, 
\begin{align} \label{005}
  \big\|[ \partial_\lambda L(\bar{\pi}_T,\lambda) ]^+ \big\| & \le  \frac{\max_{\lambda \in \Lambda_M} \|\lambda-\lambda_0\|^2 }{2r \cdot \sum_{m=0}^{T-1}  \eta_m} + G^2\frac{\sum_{m=0}^{T-1} \eta_m^2 /2}{2r \cdot\sum_{m=0}^{T-1} \eta_m} 
 + \frac{\sum_{m=0}^{T-1} \eta_m \cdot \big(L(\pi_m,\lambda_m)-L^*\big)}{2r \cdot \sum_{m=0}^{T-1} \eta_m}.
\end{align}
%It remains to upper bound the last part in the right-hand side of inequality \eqref{005}. 
Next, recall that $\bar{\pi}_T=\sum_{m=0}^{T-1} \tilde{\eta}_m \pi_m $, where
$
\tilde{\eta}_m={\eta_m}/({\sum_{m=0}^{T-1} \eta_m}),\ m=0,\ldots,T-1.
$
Since $\lambda^*$ is the optimal solution of the dual problem and $L(\pi^*,\lambda^*)\ge L(\pi^*,\bar{\lambda}_T)$, by the saddle point property, we have 
\begin{align} \label{006}
 \sum_{m=0}^{T-1} \tilde{\eta}_m \cdot \big(L(\pi_m,\lambda_m)-L^*\big) & = \sum_{m=0}^{T-1} \tilde{\eta}_m\cdot L(\pi_m,\lambda_m) -L^* \notag \\
  &\le \sum_{m=0}^{T-1} \tilde{\eta}_m\cdot L(\pi_m,\lambda_m) -L(\pi^*,\bar{\lambda}_T) \notag\\
  &{=} \sum_{m=0}^{T-1} \tilde{\eta}_m\cdot \big(L(\pi_m,\lambda_m) -L(\pi^*,{\lambda}_m)\big).
\end{align}
Similarly, under Assumption \ref{ass2}, by the second part in Lemma \ref{descent_lemma}, we have
\begin{align}  \label{007}
  \sum_{m=0}^{T-1} \tilde{\eta}_m \cdot \big(L(\pi_m,\lambda_m)-L(\pi^*,{\lambda}_m)\big) & \le \Big((1-\gamma) \cdot  \sum_{m=0}^{T-1} \eta_m \Big )^{-1}\cdot \Big (\frac{G^2}{8}\cdot  \sum_{m=0}^{T-1} \eta^2_m  + \Phi^{\pi^*}(\pi^*\| \pi_0) \Big), 
\end{align}
as the weighted KL divergence $ \Phi^{\pi^*}(\cdot||\cdot)$ is nonnegative. 

Lastly, combining inequalities \eqref{005}-\eqref{007}, we have 
\begin{align*} 
  \big\|[ \partial_\lambda L(\bar{\pi}_T,\lambda) ]^+ \big\|  \le \frac{G^2}{2r \cdot\sum_{m=0}^{T-1}\eta_m} +  \Big(\frac{1}{2}+\frac{1}{8(1-\gamma)}\Big)  G^2\frac{\sum_{m=0}^{T-1}\eta_m^2 }{2r \cdot\sum_{m=0}^{T-1} \eta_m} 
+\frac{(1-\gamma)^{-1} \Phi^{\pi^*}(\pi^*\| \pi_0)}{2r \cdot \sum_{m=0}^{T-1} \eta_m}.
\end{align*}
If we set $\eta_m={\Theta}(1/\sqrt{m})$, there exists finite constants $\kappa_1$ and $\kappa_2$ such that 
$$
\sum_{m=0}^{T-1} \eta_m \ge \kappa_1 \sqrt{T} \text{ and } \sum_{m=0}^{T-1} \eta^2_m \le \kappa_2 \log(T).
$$
Subsequently, we obtain
\begin{align*} 
  \big\|[ \partial_\lambda L(\bar{\pi}_T,\lambda) ]^+ \big\| & \le  \Big(G^2  \cdot \Big( 1+ \frac{5}{8}\kappa_2 \log(T) \Big)+ \Phi^{\pi^*}(\pi^*\| \pi_0) \Big )\frac{1}{2r(1-\gamma)\kappa_1 \sqrt{T}}.
\end{align*}
Similarly, if we set $\eta_m=\eta$ (constant step size), then  
\begin{align*} 
  \big\|[ \partial_\lambda L(\bar{\pi}_T,\lambda) ]^+ \big\| & \le  \big(G^2 +  (1-\gamma)^{-1} \cdot \Phi^{\pi^*}(\pi^*\| \pi_0) \big ) \frac{1}{2r T\eta}+  \Big(\frac{1}{2}+\frac{1}{8(1-\gamma)}\Big)\frac{G^2\eta }{2r},
 \end{align*}

{\bf We next prove the bound for $C(\bar \pi_T)-L^*$.}
We start with the upper bound. By the definition of $\bar{\pi}_T$, we have   
\begin{align} \label{008}
  C(\bar{\pi}_T)-L^*= \sum_{m=0}^{T-1}\tilde{\eta}_m \cdot \big(L(\pi_m,\lambda_m) - L^*  \big) -\sum_{m=0}^{T-1} \tilde{\eta}_m \cdot\lambda_m^\top (D(\pi_m)-q).
\end{align}
From inequalities \eqref{006} and \eqref{007}, we have 
$$
 \sum_{m=0}^{T-1}\tilde{\eta}_m \cdot \big(L(\pi_m,\lambda_m) - L^* \big) \le \Big((1-\gamma) \cdot \sum_{m=0}^{T-1} \eta_m \Big )^{-1}\cdot \Big ( \frac{G^2}{8} \sum_{m=0}^{T-1} \eta^2_m  + \Phi^{\pi^*}(\pi^*\| \pi_0) \Big).
$$
%It remains to upper bound the second part in the right-hand side of inequality \eqref{008}. 
Next, since $ D(\pi_m)-q=\partial_{\lambda}L(\pi_m,\lambda_m)$,  setting $\lambda=0$ in \eqref{002}, similarly, we obtain
$$
- \sum_{m=0}^{T-1} \tilde{\eta}_m \cdot\lambda_m^\top (D(\pi_m)-q) \le \frac{\| \lambda_0 \|^2 +G^2\cdot  \sum_{m=0}^{T-1} \eta_m^2}{2  \sum_{m=0}^{T-1} \eta_m }.
$$
Hence, if $\eta_m=\Theta(1/\sqrt{m})$, we have 
$$
C(\bar{\pi}_T)-L^* \le  \Big( \frac{5G^2}{8}\cdot \kappa_2 \cdot \log(T) +\Phi^{\pi^*}(\pi^*\| \pi_0) +\frac{\| \lambda_0 \|^2}{2} \Big)\frac{1}{(1-\gamma)\kappa_1\sqrt{T}}.
$$

For the lower bound, by the saddle point property, we have  
$$C(\bar{\pi}_T)= L(\bar{\pi}_T,\lambda^*)-(\lambda^*)^\top D(\bar{\pi}_T)\ge  L^*-(\lambda^*)^\top D(\bar{\pi}_T) . $$
Since $\lambda^* \ge 0$ and $D(\bar{\pi}_T) \le [D(\bar{\pi}_T)]^+ $,  
\begin{align*}
C(\bar{\pi}_T) - L^* & \ge  - \| \lambda^* \|  \big\|[D(\bar{\pi}_T)]^+ \big\|\\
&\ge - \| \lambda^* \|  \Big(G^2  \Big( 1+ \frac{5}{8}\kappa_2 \log(T) \Big)+ \Phi^{\pi^*}(\pi^*\| \pi_0) \Big )\frac{1}{2r(1-\gamma)\kappa_1 \sqrt{T}}.
\end{align*}
Similarly, when if $\eta_m=\eta$, we have  
\begin{align*}
C(\bar{\pi}_T)-L^*  &\le   \big(  (1-\gamma)^{-1} \Phi^{\pi^*}(\pi^*\| \pi_0) +\| \lambda_0 \|^2/2 \big)\frac{1}{T\eta} + \frac{5G^2\eta}{8(1-\gamma)},  \\
C(\bar{\pi}_T)-L^*& \ge -\| \lambda^*\|  \big(G^2 +  (1-\gamma)^{-1}  \Phi^{\pi^*}(\pi^*\| \pi_0) \big )  \frac{1}{2r T\eta} -\| \lambda^*\|  \Big(\frac{1}{2}+\frac{1}{8(1-\gamma)}\Big) \frac{G^2\eta }{2r}.
\end{align*}

%Above all, if we ignore all the constants, we have 
%\begin{align*}
% \big\| [D_{\mu_0}(\bar{\pi}_T)]^+\big\|, \big | C_{\mu_0}(\bar{\pi}_T)- L^*_{\mu_0} \big| = O\big(\log(T)/\sqrt{T}\big),
%\end{align*}  
%for $ O(1/\sqrt{m})$ stepsize and 
%\begin{align*} 
%\big\| [D_{\mu_0}(\bar{\pi}_T)]^+\big\|, \big | C_{\mu_0}(\bar{\pi}_T)- L^*_{\mu_0} \big| = O\big(1/(\eta T)+\eta \big),
%\end{align*}
%for constant stepsize $\eta$, which concludes the proof of Theorem \ref{thm1}.
\Halmos
\endproof
\end{APPENDIX}

% Appendix here
% Options are (1) APPENDIX (with or without general title) or
%             (2) APPENDICES (if it has more than one unrelated sections)
% Outcomment the appropriate case if necessary
%

% Acknowledgments here
%\ACKNOWLEDGMENT{\color{red} The authors gratefully acknowledge the existence of
%the Journal of Irreproducible Results and the support of the Society
%for the Preservation of Inane Research.}

% References here (outcomment the appropriate case)

% CASE 1: BiBTeX used to constantly update the references
%   (while the paper is being written).
%\bibliographystyle{informs2014} % outcomment this and next line in Case 1
%\bibliography{<your bib file(s)>} % if more than one, comma separated

% CASE 2: BiBTeX used to generate mypaper.bbl (to be further fine tuned)
%\input{mypaper.bbl} % outcomment this line in Case 2

%If you don't use BiBTex, you can manually itemize references as shown below.

%\bibliographystyle{nonumber}

%%%%%%%%%%%%%%%%%
\end{document}